\documentclass[english,12pt,leqno]{article}
\usepackage{tikz}
\usetikzlibrary{decorations.markings}
\usetikzlibrary{calc}
\usepackage[latin1]{inputenc}
\usepackage{amsxtra}
\usepackage{amsmath}
\usepackage{amssymb}
\usepackage{amsfonts}
\usepackage[all]{xy}
\usepackage{mathrsfs}
\usepackage{amsthm}

\usepackage{babel}

  \newenvironment{dedication}
        {\vspace{3ex}\begin{quotation}\begin{center}\begin{em}}
        {\par\end{em}\end{center}\end{quotation}\vspace{5ex}}

\newtheorem{thm}[equation]{Theorem}

\newtheorem{cor}[equation]{Corollary}
\newtheorem{lem}[equation]{Lemma}
\newtheorem{prop}[equation]{Proposition}

\newtheoremstyle{example}{\topsep}{\topsep}%
     {}
     {}
     {\bfseries}
     {.}
     {2pt}
     {\thmname{#1}\thmnumber{ #2}\thmnote{ #3}}

   \theoremstyle{example}
   
   \newtheorem{Defi}[equation]{Definition}
   
   \newtheorem{rem}[equation]{Remark}
   \newtheorem{rems}[equation]{Remarks}
   
   \newtheorem{exas}[equation]{Examples}
   \newtheorem{ex}[equation]{Example}

\newtheorem{conv}[equation]{Convention}
\newtheorem{CYP}[equation]{2-Calabi-Yau principle}

\newcommand{\Par}{\vspace{5pt}\noindent\textbf{ 
\refstepcounter{equation}
(\theequation)}
  \textbf}
  
  \newcounter{sparagraphcounter}[section]

\def\sparagraph#1{\refstepcounter{sparagraphcounter}
\vskip 1em%
\noindent {\bf \Alph{sparagraphcounter}. #1}%
}

\newtheoremstyle{example}{\topsep}{\topsep}%
     {}
     {}
     {\bfseries}
     {.}
     {2pt}
     {\thmname{#1}\thmnumber{ #2}\thmnote{ #3}}

   \numberwithin{equation}{section}

\def\eps{{\varepsilon}}

\def\vpi{\varpi}

\def\AAA{\mathbb{A}}

\def\CC{\mathbb{C}}
\def\DD{\mathbb{D}}
\def\FF{\mathbb{F}}
\def\GG{\mathbb{G}}
\def\LL{\mathbb{L}}
\def\PP{\mathbb{P}}
\def\RR{\mathbb{R}}
\def\ZZ{\mathbb{Z}}

\def\HH{\mathbb{H}}

\def\gen{\mathfrak{g}}

\def\Ien{\mathfrak{I}}
\def\Jen{\mathfrak{J}}

\def\men{\mathfrak{m}}
\def\Men{\mathfrak{M}}

   \def\Yen{\mathfrak{Y}}

\def\Ac{\mathcal{A}}
\def\Bc{\mathcal{B}}
\def\Cc{\mathcal{C}}

\def\Dc{\mathcal{D}}
\def\Ec{\mathcal{E}}
\def\Fc{\mathcal{F}}
\def\Gc{\mathcal{G}}

\def\Lc{\mathcal{L}}
\def\Mc{\mathcal{M}}
\def\Nc{\mathcal{N}}

\def\Oc{\mathcal{O}}

\def\kb{\mathbf{k}}

\def\1{{\bf 1}}

\def\lsb{[\hskip -.04cm [}

\def\be{\begin{equation}}

\def\dim{{\rm{dim}}}
\def\Dbc{D^b_{\on{constr}}}
\def\DSh{D^b\Sh }
\def\Dbcc{D^b_{\on{cc}}}

\def\ee{\end{equation}}
\def\Ext{\on{Ext}}

\def\FT{\on{FT}}

\def\Hom{\on{Hom}}
\def\Homu{\ul\Hom}

\def\k {\mathbf k}

\def\lla{\longleftarrow}

\def\lra{\longrightarrow}
\def\lsb{[\hskip -.04cm [}
\def\lrb{(\hskip -.07cm (}

\def\mhom{\mu\ul{\on{Hom}}}
\def\Mhom{\Mc\ul\Hom}
\def\Mod{\on{Mod}}

\def\Ob{\on{Ob}}
\def\ol{\overline}
\def\on{\operatorname}

\def\Perv{\on{Perv}}

\def\q{{\mathbf q}}
\def\qq{ \slash \hskip -.12cm \slash }
\def\qqq{\slash  \hskip -.12cm \slash \hskip -.12cm \slash }
\def\qpm{{\q^{\pm 1}}}

\def\rsb{]\hskip -.04cm ]}
\def\rrb{)\hskip -.07cm )}

\def\Sh{\on{Sh}}

\def\Spec{{\on{Spec}}}

\def\ul{\underline}

\def\Vect{\on{Vect}_\k}

\def\wt{\widetilde}

\newcommand{\F}{{\mathcal F}}

\title{Microlocal sheaves and quiver varieties}

\author{{Roman Bezrukavnikov} \and {Mikhail Kapranov} }

\begin{document}

\maketitle

\begin{dedication} \`A Vadim Schechtman pour son 60-\`eme anniversaire 
\end{dedication}

 


\tableofcontents


 \addtocounter{section}{-1}

\section{Introduction.}

The goal of this paper is to relate  two classes of symplectic manifolds of great importance in Representation Theory
and to put them into a common framework.

\vskip .2cm

\Par{ Moduli  of local systems on Riemann surfaces.} 
First, let $X$ be a compact oriented  $C^\infty$ surface  and $G$ be a reductive algeraic group. 
The moduli space $\on{LS}_G(X)$ of $G$-local systems on $X$ is naturally a symplectic manifold \cite{goldman},
with the symplectic structure given by the cohomological pairing. 
As shown by Atiyah-Bott, $\on{LS}_G(X)$ can be obtained as the Hamiltonian reduction of an infinite-dimensional
flat symplectic space formed by all $G$-connections, with the Lie algebra-valued  moment map given by the curvature. 
Alternatively,  $\on{LS}_G(X)$
  can be  obtained as a Hamiltonian reduction of a finite-dimensional 
   symplectic space but at the price
of passing to the { multiplicative theory}: replacing the Lie algebra-valued moment map by a group-valued one \cite{AMM}. 

The variety  $\on{LS}_G(X)$ and its versions associated to surfaces with punctures, marked points etc.
form fundamental examples of cluster varieties \cite{FG}, and their quantization is interesting from many points of view. 
We will be particularly interested in the case $G=GL_n$, in which case local systems form an abelian category.

\vskip .2cm

\Par{ Quiver varieties.} The second class is formed by the  Nakajima quiver varieties \cite{N}. 
Given a finite oriented graph $Q$, the corresponding quiver varieties can be seen as symplectic reductions of
the cotangent bundles to the moduli spaces of representations of $Q$ with various dimension vectors. 
Passing to the cotangent bundle has the effect of  ``doubling the quiver": 
introducing, for each arrow $\xymatrix{ i \ar[r] &j}$ of $Q$, a  new arrow $\xymatrix {i & \ar[l] j} $ in the opposite direction.

 Interestingly, one also has 
 the  ``multiplicative" versions of quiver varieties defined  by Crawley-Boevey
and Shaw \cite{CBS} and Yamakawa \cite{yamakawa}. They can   be constructed by performing the 
Hamiltonian  reduction  but using the group-valued moment map. 
It is these multiplicative versions that we will consider in this paper.

\vskip .2cm

\Par { Relation to perverse sheaves.} 
It turns out that both these classes  can be put under the same umbrella of {\em varieties arising from classification of perverse sheaves}.

From the early days of the theory \cite{BBD}, a lot of effort has been spent on finding   descriptions of various categories of perverse sheaves
as representation categories of some explicit quivers with relations. In all  of these cases, the quivers have the following remarkable
property: {\em their arrows come in pairs of opposites}
$ \xymatrix{
 i
   \ar@<.4ex>[r]&  \ar@<.4ex>[l]
  j
  }
  $.
  This reflects the fact that any category of perverse sheaves has a perfect duality (Verdier duality). The  diagram (representation  of the quiver)
    corresponding to
  the dual perverse sheaf $\Fc^\bigstar$ is obtained  from the diagram corresponding to $\Fc$
  by dualizing both  the spaces and (up to a minor twist, cf. \cite[(II.3.4)]{Ma}) the arrows, thus interchanging the elements of each  pair
  of opposites. 
  We see therefore a conceptual reason for a possible relationship between perverse sheaves and quiver varieties. 
  
  The relation between perverse sheaves and  $\on{LS}_{GL_n}(X)$  is even more immediate: local systems are nothing but
  perverse sheaves without singularities, so ``moduli spaces of perverse sheaves" are natural objects to look at. 
  
  \vskip .2cm
  
  \Par { Microlocal sheaves.} 
  However, to make the above relations precise, we need to use a generalization of perverse sheaves: {\em microlocal sheaves}.
  These objects can be thought as modules over a (deformation) quantization of  a symplectic manifold $S$ supported in a given
  Lagrangian subvariety $X$, see \cite{KS-DQ}. The case $S=T^*M$ being the cotangent bundle to a manifold $M$ and $X$ being
  conic, corresponds to the usual theory of holonomic $\Dc$-modules and perverse sheaves. However, for our applications it is
  important to consider the case when $X$ is compact. 
  
  In this paper we need only the simplest case when $X$ is an algebraic curve over $\CC$ which is allowed to have
  nodal singularities. In this case microlocal sheaves can be defined in a very elementary way as perverse sheaves on
  the normalization satisfying a Fourier transform condition near each self-intersection point. 
  The relation with quiver varieties appears when we take $X$ to be a union of projective lines whose intersection graph 
  is our ``quiver" $Q$ (with orientation ignored). 
  
  If we consider only ``smooth" microlocal sheaves (no singularities other than the nodes), we get a natural analog of the concept
  of a local system for nodal curves.  In particular, for a compact $X$ 
we consider such microlocal sheaves as objects of a triangulated category $D\Mc(X,\emptyset)$
of {\em microlocal complexes}, and we show in Thm. \ref{thm:CY} that it has the {\em 2-Calabi-Yau
property},  extending  the Poincar\'e duality for local systems:
\[
R\Hom (\Fc, \Gc)^* \,\,\simeq \,\, R\Hom(\Gc, \Fc)[2].
\]
This gives an intrinsic reason to expect that the ``moduli spaces"  parametrizing
microlocal sheaves or complexes, are symplectic, 
 in complete analogy with Goldman's picture \cite{goldman} for local systems. 
 We discuss the related   issues in \S \ref{sec:preprogen}D and give a more direct
 construction of such spaces in \S \ref{sec:framed} by using quasi-Hamiltonian reduction.

  \vskip .2cm
  
  \Par{Relation to earlier work.} An earlier attempt to relate (multiplicative)  quiver varieties and D-module type objects (i.e., to invoke the Riemann-Hilbert correspondence) was made by D. Yamakawa \cite{yamakawa}. Although his construction is  quite different from ours and is  only 
  applicable to quivers of a particular shape, it was one of the starting points of our inverstigation. 
  
  \vskip .2cm
  
  More recently, a Riemann-Hilbert type interpretation of multiplicative 
  preprojective algebras was given by W.~Crawley-Boevey \cite{CB}.
  His setup is in fact quite close to ours (although we learned of his paper only after
  most of  our constructions have been formulated). In particular, the datum of  a
  ``Riemann surface quiver with non-interfering arrows", a central concept of \cite{CB}, is equivalent to the datum of a nodal curve $X$: the normalization
  $\wt X$ is then the corresponding Riemann surface, and the pairwise identifications of the points of $\wt X$  needed to get $X$, form a Riemann surface quiver. From our point of view, the construction of \cite{CB} can be seen as leading to
    an explicit description, in terms of D-module type data,
   of ``smooth" microlocal sheaves on  a nodal curve, see Theorem \ref{thm:micro-dr}. 
   
   \vskip .2cm

 Considering a nodal curve $X$ as the basic object,
  has  the advantage of putting the situation,
  at least heuristically,  into the general framework of
  deformation quantization (DQ-)modules. In particular, one can consider for $X$ a
 projective  curve with more complicated singularities,  realized as a (necessarily Lagrangian) subvariety
  in a holomorphic symplectic surface. The general theory of \cite{KS-DQ}
  suggests that moduli spaces of ``smooth"  microlocal sheaves in this
situation will produce interesting symplectic varieties.   
  Further, passing to higher-dimensional projective 
  singular Lagrangian
  varieties $X$, one expects to get
  shifted symplectic varieties,  as suggested  by the Calabi-Yau property
  of DQ-modules \cite[Cor. 6.2.5]{KS-DQ} and  the general theory of \cite{kontsevich-soibelman}
 and   \cite{PTVV}.

\vskip .2cm

\Par { Acknowledgements.} We are grateful to A. Alekseev, Y. Brunebarbe, T. Dyckerhoff,
V. Ginzburg, P. Schapira, 
Y. Soibelman and G. Williamson for useful discussions and correspondence.
  The work of M.K.  was supported by the  World Premier International Research Center Initiative (WPI), MEXT, Japan
and     by the Max-Planck Institute  f\"ur Mathematik, Bonn, Germany.

\vskip .3cm

\vfill\eject

\vfill\eject

\section{Microlocal sheaves on nodal curves}\label{sec:microsheaves}

\sparagraph{Topological definitions.} Let $X$ be a nodal 
curve 
  over $\CC$, i.e., an algebraic, quasi-projective curve  whose only
  singularities are transversal self-intersection points 
  (also known as {\em nodes}, or ordinary double points).

  For a node $x\in X$ we denote   
two ``branches" of $X$ near $x$ (defined up to permutation) by $B'$ and $B''$. 
More precisely, we think of $B'$ and $B''$ as small disks meeting at $x$. 
Alternatively, let $\vpi: \wt X\to X$ be the normalization of $X$. Then 
$\vpi^{-1}(x)=\{x', x''\}$ consists of two points, and we define 
$\wt {B}', \wt {B}''$ as the neighborhoods of $x'$ and $x''$ in $\wt X$. 
We can then identify canonically $B' = \wt{B}'$, $B''=\wt{B}''$. 
We note that the Zariski tangent space to $X$ at a node $x$ is 2-dimensional:
\[
T _xX \,\,=\,\, T_x B' \oplus T_x B''. 
\]

 \begin{Defi}
 A {\em duality structure} on $X$ is a datum, for each node $x$, of a symplectic structure $\omega_x$ on the 2-dimensional vector space $T_xX$. 
 \end{Defi}

Alternatively, a duality structure at a node $x$ can be considered as a
datum of isomorphisms
\[
\eps'_x: T_x B' \to T^*_x B'', \quad \eps_x'': T_x B'' \to T^*_x B'
\]
such that $(\eps''_x)^* = - \eps'_x$. 

\begin{ex}\label{ex:SS} (a) 
Suppose $X$ embedded into a holomorphic symplectic surface $(S,\omega)$. Then the restrictions
of $\omega$ to all the nodes of $X$ give a duality structure on $X$. 

Note that any duality structure on $X$ can be obtained in this way.
 Indeed, we first consider a neighborhood $\wt S$ of the zero section in the cotangent 
 bundle $T^*\wt X$. Then  for any node $x\in X$ with $\vpi^{-1}(x)=\{x', x''\}$, 
 we identify the neighborhoods $U'$ of $x'$ and $U''$ of 
  $x''$ in $\wt S$ by an appropriate symplectomorphism
   so that the intersection of $U'$ with the  zero section of $T^*\wt S$ becomes
   identified with the intersection of $U''$ with the fiber of $T^*\wt S$ over $x''$
   and vice versa. 
   
   \vskip .1cm
   
   (b) Situations when $X$ is naturally
   embedded into an {\em algebraic} symplectic surface $S$, provide a 
   richer structure. The best known examples are provided by  $S$ being the  minimal resolution
   of a Kleinian singularity $\CC^2/G$, where $G$
   is a finite subgroup in $SL_2(\CC)$.
   In this case $X$ is a union of projective lines,
   with the intersection graph being a Dynkin diagram of type ADE.

\end{ex}

Let $X$ be a nodal curve with a duality structure. For each node $x\in X$
 we can identify $B'$ and $B''$ with open disks in $T_xB'$ and $T_x B''$
 or, equivalently, in $T_x\wt{B}'$ and $T_x {\wt B}''$ respectively. Such identifications
 are unique up to contractible spaces of choices.

  Let $D^b(\wt {B}', x')$ be the
  full subcategory in $\Dbc(\wt {B}')$ formed by complexes whose cohomology sheaves
  are locally constant outside $x'$, and similarly for  $D^b(\wt {B}'', x'')$.
 Let $\Perv(\wt{B}', x')\subset D^b(\wt {B}', x')$
and $\Perv(\wt{B}'', x'')\subset D^b(\wt {B}'', x'')$
be the full (abelian)  subcategories formed by  perverse sheaves.

  The above identifications with the disks in the tangent spaces together with the
  isomorphisms $\eps', \eps''$ give rise to geometric Fourier(-Sato) transforms which are
  equivalences of pre-triangulated categories
   \be\label{eq:FT-1}
     \xymatrix{
   D^b(\wt{B}', x')  
   \ar@<.7ex>[rr]^{\on{FT}'}&&  \ar@<.7ex>[ll]^{\on{FT}''}
    D^b(\wt{B}'', x''),
  }
  \ee
which are canonically inverse to each other and restrict to equivalence of abelian categories
 \be\label{eq:FT-2}
     \xymatrix{
   \Perv(\wt{B}', x')  
   \ar@<.7ex>[rr]^{\on{FT}'}&&  \ar@<.7ex>[ll]^{\on{FT}''}
    \Perv(\wt{B}'', x''). 
  }
  \ee

 \begin{rem}
 The fact that $\on{FT}'$ and $\on{FT}''$ are precisely inverse to each other, comes from the requirement that 
 $\eps'_x$ and $\eps''_x$ are  the negatives of the transposes of each other, rather than exact transposes.
  We recall that the ``standard" Fourier-Sato transform for a $\CC$-vector space
 $E$ is an equivalence (\cite{KS}, Ch. 3)
 \[
 \on{FT}_E: D^b_{\on{mon}}(E) \to D^b_{\on{mon}}(E^*)
 \]
 ($D^b_{\on{mon}}$ means the derived category of $\CC$-monodromic constructible complexes). In this setting
 $\on{FT}_{E^*}$ is not canonically inverse to $\on{FT}_E$: the composition
 $\on{FT}_{E^*}\circ \on{FT}_E$ is canonically identified with $(-1)^*$, the pullback
 with respect to the antipodal transformation $(-1): E\to E$.  
 \end{rem}

\begin{Defi}\label{def:micro-complex}
A {\em microlocal complex} $\Fc$ on $X$ is a datum consisting of:
\begin{enumerate}
\item [(1)] A $\CC$-constructible complex $\wt\Fc$ on $\wt X$.

\item[(2)] For each node $x\in X$,   quasi-isomorphisms of constructible complexes 
\[
\alpha': \wt\Fc|_{\wt{B}'} \lra \on{FT}''\left( \wt\Fc|_{\wt{B}''}\right),\quad 
\alpha'': \wt\Fc|_{\wt{B}''} \lra \on{FT}'\left( \wt\Fc|_{\wt{B}'}\right),
\]
 inverse to each other. 
\end{enumerate}
A {\em microlocal sheaf} on $X$ is a microlocal complex $\Fc$ such that $\wt\Fc$
is a perverse sheaf on $\wt X$. 
\end{Defi}

A {\em morphism} of microlocal complexes (resp. microlocal sheaves)
 $\Fc\to\Gc$ is a morphism of constructible complexes (resp. perverse sheaves)
$\wt\Fc\to\wt\Gc$ on $\wt X$ compatible with the identifications $\alpha', \alpha''$. 
In this way we obtain a pre-triangulated category
  $D\Mc(X)$  formed by microlocal complexes on $X$ and an abelian subcategory
  $\Mc(X)$ formed by microlocal sheaves. 
  
  For a finite subset of smooth points $A\subset X_{\on{sm}}$ we denote by $D\Mc(X, A)\subset D\Mc(X)$
  the full subcategory formed by microlocal complexes $\Fc$ such that $\wt\Fc$ is
  smooth (i.e.,  each cohomology sheaf of it is a local system) 
  outside of $\vpi^{-1}(A)$. Let $\Mc(X,A)$ be the intersection of $\Mc(X)$
  with $D\Mc(X,A)$. 
  
   \begin{rems}\label{rems:DQ}
    (a) 
  Suppose $\k=\CC \lrb h\rrb$ is the field of Laurent series in one variable $h$ with
  complex coefficients. 
   Assume that $X$ is embedded into a symplectic surface $(S,\omega)$, as in
  Example \ref{ex:SS}. As shown in \cite{KS-DQ}, $S$ admits a 
  {\em deformation quantization
  algebroid} $\Ac_S$, which locally can be viewed as a sheaf of $\CC \lsb h\rsb$-algebras
  whose reduction modulo $h$ is identified with $\Oc_S$ and whose first order
  commutators are given by the Poisson bracket of $\omega$. One also has the
  $h$-localized algebroid 
  $\Ac_S^{\on{loc}}=\Ac_S \otimes_{\CC \lsb h\rsb} \CC\lrb h
  \rrb$. 
  
  \vskip .2cm
  
  The category $D\Mc(X, \emptyset)$ can be compared with the category 
  $D^b_{\on{gd}, X}(\Ac_S^{\on {loc}})$ of complexes of $\Ac_X^{\on{loc}}$-modules whose
  cohomology modules are coherent, algebraically good \cite[2.7.2]{KS-DQ}
  modules
   supported on $X$. More precisely, each smooth (not necessarily closed)
   Lagrangian $\CC$-submanifold (i.e., a smooth complex curve) $\Lambda\subset S$,
   gives a simple holonomic  $\Ac_X^{\on{loc}}$-module $\Oc_\Lambda$,
   and we have the  ``$\Lambda$-Riemann-Hilbert functor" 
   \[
   R\Homu_{\Ac_S^{\on {loc}}}(-, \Oc_\Lambda): D^b_{\on{gd}, X}(\Ac_S^{\on {loc}})\lra \Dbc(\Lambda).
   \]
   Taking for $\Lambda$ various smooth branches of $X$, we associate to an object $\Nc$
   of $D^b_{\on{gd}, X}(\Ac_S^{\on {loc}})$ a constructible complex  $\wt\Fc$ on
   $\wt X$. If   $\Nc$ is a single module
   in degree 0, then $\wt\Fc$ is a perverse sheaf.
    When two branches meet at a point (node $x$ of $X$),
    the corresponding Riemann-Hilbert functors 
   are, near $x$,  related to each other by the Fourier transform, thus leading to Definition 
   \ref {def:micro-complex}. 
  
   \vskip .2cm
   
    (b) A particularly interesting algebraic case is provided by $S$ being the
   minimal resolution of a Kleinian singularity, see Example 
    \ref{ex:SS}(b). In this case quantizations of $S$
exist algebraically in finite form
 (not just over power series in $h$), see \cite{Bo}.
It is therefore interesting to compare their modules with microlocal sheaves
on Dynkin chains of $\PP^1$'s.   
  \end{rems} 

  Let $X$ be a nodal curve with  duality structure and $A\subset X_{\on{sm}}$ a finite
  subset of smooth points. Let us form a new, noncompact nodal curve
  \[
  X_A \,\,=\,\, X \,\cup \, \bigcup_{a\in A} T^*_aX
  \]
  by attaching each cotangent line $T^*_aX$ to $X$ at  the point $a$
  which becomes a new node. The symplectic
  structure on $T^* X_{\on{sm}}$ gives a duality structure at each new node.
  
  \begin{prop}\label{prop:A>empty}
  We have canonical equivalences\[
  D\Mc(X,A) \,\,\simeq \,\, D\Mc(X_A,\emptyset), \quad  \Mc(X,A) \,\,\simeq \,\, \Mc(X_A,\emptyset). 
  \]
  \end{prop}
  
  \noindent{\sl Proof:} We identify the normalization of $X_A$ as
  \[
  \wt X_A\,\,=\,\, \wt X \, \sqcup \,\bigsqcup_{a\in A} T^*_aX. 
  \]
  To each microlocal complex $\Fc$ on $X$ we associate a microlocal complex $\Fc_A$ on $X_A$ given by
  \[
  \wt\Fc_A|_{\wt X}=\wt\Fc,\quad \wt\Fc_A|_{T^*_aX} = \mu_a(\Fc), 
  \]
   where $\mu_a(\Fc)$ is the microlocalization of $\Fc$ at $a$, 
   i.e., the Fourier transform of the specialization of $\Fc$ at $a$ \cite{KS}.
  The definition gives the Fourier transform identifications for $\wt\Fc_A$. 
    This defines the desired equivalence.
   \qed
  
  \sparagraph{The Calabi-Yau property.} Important for us will be the following.

  \begin{thm}\label{thm:CY}
  Let $X$ be a compact nodal curve over $\CC$ equipped with a duality
  structure. Then $D\Mc(X,\emptyset)$ is a Calabi-Yau dg-category of dimension 2.
  In other words, for any $\Fc, \Gc\in D\Mc(X,\emptyset)$ we have a canonical quasi-isomorphism of complexes of $\k$-vector spaces
  \[
  R\Hom(\Fc, \Gc)^* \simeq R\Hom(\Gc, \Fc)[2].
  \]
  \end{thm}
  
   \begin{ex}
  For   $X$ smooth, the category $\Mc(X,\emptyset)$ consists of local
  systems on $X$, and $D\Mc(X)$ consists of complexes with locally constant cohomology. 
  Theorem \ref{thm:CY} in this case reduces to  the Poincar\'e duality for local systems on
  a compact oriented topological surface.   \end{ex}
  
  \begin{rem} 
  Consider the situation of Remark \ref{rems:DQ}(a). 
 For a compact symplectic manifold $S$ of any dimension $d$, Corollary
   6.2.5 of \cite{KS-DQ} gives that $D^b_{\on{gd}}(\Ac_S^{\on {loc}})$,
   the category of all complexes of  $\Ac_S^{\on {loc}}$-modules with
   coherent and  algebraically good
   cohomology, is a  Calabi-Yau category over $\CC\lrb h\rrb$ of dimension $d$.
   This result can be seen as a noncommutative lifting of the classical Serre duality
   for coherent $\Oc_S$-modules. 
  
    \vskip .2cm
   
    If $S$ is non-compact, then restricting the support to a given
   compact subvariety $X$ allows one to preserve the duality,
   cf. \cite[Cor. 3.3.4]{KS-DQ}. In particular, when $S$ is a symplectic surface,
   and $X\subset S$ is a compact nodal curve, $D^b_{\on{gd}, X}(\Ac_S^{\on {loc}})$
   is a Calabi-Yau category over $\CC\lrb h\rrb$ of dimension 2. Our Theorem
   \ref{thm:CY} can be seen as a topological analog of this fact. 
     \end{rem}
     
     \noindent {\bf Proof of Theorem \ref{thm:CY}.} Let $\Fc, \Gc\in
 D\Mc(X,\emptyset)$. For any open set $U\subset X$ (in the classical topology) we
  have the complex of vector spaces
 \[
 R\Hom_{D\Mc (U, \emptyset)}(\Fc|_U, \Gc|_U)
 \,\,\,\in \,\,\, D^b\Vect. 
 \]
Taken for all $U$, these complexes can be thought as forming a
complex of sheaves which we denote 
\[
\Mhom (\Fc, \Gc) \,\,\,\in \,\,\, \Dbc (X),
\]
so that, in a standard way, we have
\[
R\Hom_{D\Mc(X,\emptyset)}(\Fc, \Gc) \,\,=\,\, R\Gamma(X, \Mhom(\Fc, \Gc)). 
\]
Our statement will follow from the Poincar\'e-Verdier duality
on the compact space $X$, if we establish the following.

\begin{prop}\label{eq:Mhom}
For any nodal curve $X$  (compact or not) with duality structure and any microlocal complexes $\Fc,\Gc
\in D\Mc(X,\emptyset)$ we have a canonical identification
\[
\DD_X \Mhom(\Fc, \Gc) \,\,\simeq \,\, \Mhom(\Gc, \Fc)[2]. 
\]
\end{prop}

To prove the proposition, we compare the bifunctor $\Mhom$ with the 
microlocal Hom bifunctor  of \cite{KS} which we recall. 

Let $M$ be a smooth manifold and $\pi: T^*M\to M$ be its cotangent bundle.
For any two complexes of sheaves $F, G$ on $M$. Kashiwara and
Schapira \cite{KS} defined a complex of sheaves
\[
\mhom(F, G) \,\,\in \,\, D^b\Sh_{T^*M}
\]
so that
\[
\begin{gathered}
R\Homu(F, G) \,\,=\,\, R\pi_*\bigl( \mhom(F, G)\bigr), 
\\
 R\Hom_{D^b\Sh_M}(F, G) \,\,=\,\, R\Gamma\bigl(T^*M, 
\mhom(F, G)\bigr).
\end{gathered}
\]

\begin{lem}\label{lem:mhom}
Assume that $M$ is a complex manifold and
$F,G\in\Dbc(M)$. Then we have a canonical identification
\[
\DD_{T^*M} \bigl( \mhom(F,G)\bigr) \,\,\simeq  \,\, \mhom(G,F) . 
\]
\end{lem}

\noindent {\sl Proof:} This is a particular case of Proposition 8.4.14(ii) of \cite{KS}. \qed

\vskip .2cm

We now deduce Proposition \ref{eq:Mhom} from Lemma \ref{lem:mhom}. 

\begin{Defi} Call a subset $Z\subset X$
{\em unibranched}, if
   $Z$ is the image, under the normalization map $\varpi: \wt X\to X$, of an open (in the classical topology)
subset $\wt Z$ such that the restriction $\varpi|_{\wt Z}: \wt Z\to Z$ is a bijection. 
\end{Defi}

Note that a unibranched subset $Z$ is a complex analytic curve which may not be open in $X$, if it passes through some nodes
(in which case it contains only one branch near each node it passes through). 
For a microlocal complex $\Fc$ on $X$ and a unibranched $Z\subset X$ we have a well-defined constructible
complex 
\[
\Fc||_Z\ \,\, := \,\,  (\varpi|_{\wt Z})_*  \wt \Fc \,\, \in\,\,  \Dbc(Z). 
\]
 Assume that $X$ is embedded into a symplectic surface $S$ and let $U$ be a neighborhood
  of $Z$ in $S$. Then we can make the following identifications:
  
  \begin{itemize}
  \item[(1)] $U$  can be identified with a neighborhood of $Z$ in $T^*Z$ so that $Z$ becomes identified
  with the zero section $T^*_ZZ$. 
  
  \item[(2)] 
   If we denote the nodes of $X$ contained in $Z$, by $x_i, i\in I$,  then $U\cap Z$
can be identified with  the union of   $T^*_Z Z$ and of some neighborhoods of 0
in the fibers $T^*_{x_i} Z$. 

\item[(3)] 
Let  $\Fc, \Gc$ be two microlocal complexes on $X$. Then, under the above identifications,
we have an isomorphism
\[
\Mhom(\Fc, \Gc)|_{U\cap Z}\,\,\simeq \,\, \mhom(\Fc||_Z, \Gc||_Z)|_{U\cap Z}. 
\]
 \end{itemize}
 
 Further, because of the Fourier transform identifications in the definition of a microlocal complex, the identifications
 in (3) are compatible for different unibranched sets passing through a given node. 
 Therefore the identifications (3) allow us to glue the identifications of Lemma \ref{lem:mhom}
 to a canonical identification as in 
 Proposition \ref{eq:Mhom}. This proposition and  Theorem \ref{thm:CY} 
 are now proved.

 \vfill\eject
 
 \section {Microlocal sheaves:   de Rham description.}\label{sec:dr}

 We now give a $\Dc$-module type description of microlocal sheaves,
 relating our approach with that of \cite{CB}. 
 
 \sparagraph{Formulations.}
 Let $X$ be a nodal curve with the set of nodes $D$ and its preimage $\wt D=\varpi^{-1}(D)\subset \wt X$.
  By an {\em orientation} of $X$ we mean a choice, for each node $x$, of the order $(x'<x'')$ on the two
  element set of preimages $\varpi^{-1}(x) = \{x', x''\}$. 
  
  We denote by
  \[
  \Re^{-1}[0,1) \,\,=\,\, [0,1) + i\RR \,\,\subset \,\, \CC
  \]
 the standard fundamental domain for $\CC/\ZZ$.
 
 Let $Y$ be a smooth algebraic curve over $\CC$ (not necessarily compact) and $Z\subset Y$ a finite subset.
 We recall, see, e.g. \cite{Ma}, the concept of a {\em logarithmic connection} (along $Z$) on an algebraic vector bundle $\Ec$
 on $Y$. Such a connection $\nabla$ can be viewed as an algebraic differential operator
  $\nabla: \Ec\to \Ec\otimes\Omega^1_Y(\log Z)$. It has a well-defined {\em residue} $\on{Res}_z(\nabla)\in \on{End}(\Ec_z)$
  at each $z\in Z$. 
  For a noncompact $Y$ there is a concept of a {\em regular} logarithmic connection (having regular singularities
  at the infinity of $Y$). 
  
  \begin{Defi}
  Let $X$ be a nodal curve over $\CC$, not necessarily compact, with orientation. A {\em de Rham microlocal sheaf}
  (without singularities)
  on $X$ is a datum of:
  \begin{enumerate}
  \item[(1)] A vector bundle $\Ec$ on $\wt X$, together with a regular logarithmic connection $\nabla$ along $\wt D$.

  \item[(2)] For each node $x\in D$ with preimages $x', x''\in\wt D$ (order given by the orientation),
  two linear operators 
  \[
   \xymatrix{
   \Ec_{x'}
   \ar@<.7ex>[r]^{u_x}&  \ar@<.7ex>[l]^{v_x}
 \Ec_{x''}
  }
  \]
  such that:
  
  \item[(3)] $\on{Res}_{x'}(\nabla) = v_x u_x, \quad \on{Res}_{x''}(\nabla) = -u_x v_x$; 
  
  \item[(4)] All eigenvalues of $v_xu_x$ and $-u_x v_x$ lie in $\Re^{-1}[0,1)$. 
   \end{enumerate}
   The category of de Rham microlocal sheaves on $X$ without singularities will be denoted by $\Mc_{\on{dR}}(X,\emptyset)$. 
  \end{Defi}
  
  \begin{rems}
  A de Rham microlocal sheaf is a particular case $(\lambda=0)$ of a $\lambda$-connection of \cite{CB},
  but with additional restriction (4). 
  \end{rems}
  
  \begin{thm}\label{thm:micro-dr}
  Take the base field $\k=\CC$. 
  Assume that $X$ is equipped with both an orientation and a duality structure. Then $\Mc_{\on{dR}}(X,\emptyset)$
  is equivalent to $\Mc(X,\emptyset)$. 
  \end{thm}

  \sparagraph{ Riemann-Hilbert correspondence.}
  In order to prove Theorem \ref{thm:micro-dr}, 
 we recall  two classical results about the Riemann-Hilbert correspondence.
  
    First, let $Y$ be a smooth curve over $\CC$ and $Z\subset Y$ a finite subset. 
A  regular logarithmic connection $\nabla: \Ec\to \Ec\otimes\Omega^1_Y(\log Z)$ will be called {\em canonical},
 all eigenvalues of all $\on{Res}_z(\nabla)$, $z\in Z$, lie in $\Re^{-1}[0,1)$. In this case $(E,\nabla)$
 is obtained by the {\em Deligne canonical extension} from its restriction to $Y-Z$, see \cite{Ma}. 
We denote by $\on{Conn}_{\on{can}}^{\on{reg}}(Y,Z)$ the category of vector bundles with 
  regular canonical connections.
  
  \begin{prop}\label{prop:can}
 The category $\on{Conn}_{\on{can}}^{\on{reg}}(Y,Z)$ is equivalence to $\on{LS}(Y-Z)$, the category
 of local systems on $Y-Z$. The equivalence is obtained by restricting $(\Ec,\nabla)$ to $Y-Z$ and taking
 the sheaf of covariantly constant sections. \qed
    \end{prop}
    
    \begin{prop}\label{prop:I}
    \cite{Ka}\cite[(II.2.1)]{Ma}
    Let $\Ien$ be the category of diagrams of finite-dimensional $\CC$-vector spaces
      \[
 H= \bigl\{ \xymatrix{
   E
   \ar@<.7ex>[r]^{u}&  \ar@<.7ex>[l]^{v}
F
  }\bigr\}
  \]
  s.t. all eigenvalues of $uv$ and $vu$ lie in $\Re^{-1}[0,1)$. Then $\Ien$ is equivalent to $\Perv(\CC,0)$. The equivalence
  takes an object $ H\in\Ien$ to the $\Dc_\CC$-module $M_H$ with the space of generators $E\oplus F$ and relations
  \[
  \begin{gathered}
  x\cdot f = v(f), \,\,\, f\in F,
  \\
  {d\over dx} \cdot e = u(e), \,\, e\in E, 
  \end{gathered}
  \]
 and then to the de Rham complex of $M_H$. \qed
    \end{prop}

\sparagraph{Fourier transform and RH.}
Recall \cite{Ma} that the Fourier-Sato transform on $\Perv(\CC,0)$
corresponds, at the $\Dc$-module level, to passing from 
passing from the generators $x, {d\over dx}$ of the Weyl algebra of  
differential
operators to new generators 
\[
p = -{d\over dx}, \,\,\, {d\over dp} = x, 
\]
so that
\[
\left[ {d\over dp}, \, p\right] \,\,=\,\, \left[ {d\over dx}, \, x\right] 
\,\,=\,\, 1. 
\]
This implies:

\begin{cor}\label{cor:FT-I}
     The effect of the Fourier-Sato transform on $\Ien$ is the functor
    \[
\FT_\Ien: \,\,\,  H =   \bigl\{  \xymatrix{
   E
   \ar@<.7ex>[r]^{u}&  \ar@<.7ex>[l]^{v} F
   }  \bigr\}
   \quad \buildrel{}\over\longmapsto \quad 
  \hat H =  \bigl\{\xymatrix{
   F
   \ar@<.7ex>[r]^{v}&  \ar@<.7ex>[l]^{-u}
E
  }\bigr\}. \quad\quad\quad \qed
  \]
\end{cor}  
  
  Therefore we can reformulate Proposition \ref{prop:I} in a more ``microlocal" form
  
  \begin{prop}
  Let $C=\{ x p=0\} \subset\CC^2$ be the coordinate cross with the orientation defined by putting the
  $x$-branch before the $p$-branch. Then $\Mc_{\on{dR}}(C,\emptyset)$ is equivalent to  $\Perv(\CC,0) \simeq \Mc(C,\emptyset)$. 
    \end{prop}
    
   \noindent {\sl Proof:} For a diagram $H\in\Ien$, the $\Dc_\CC$-module $M_H$ becomes $\Oc$-coherent
   on $\CC-\{0\}$, and is identified with the following bundle with connection:
   \[
   \Ec_{H}^0 \,\,=\,\,\biggl( E\otimes\Oc_{\CC-\{0\}}, \,\nabla = d- (vu) {dx\over x}\biggr). 
   \] 
   Therefore the Deligne canonical extension of $\Ec^0_H$ to $\CC$ is the logarithmic connection
   \[
   \Ec_{H} \,\,=\,\,\biggl( E\otimes\Oc_{\CC}, \,\nabla = d- (vu) {dx\over x}\biggr). 
   \]  
   Similarly for the Fourier transformed diagram $\hat H$ which gives a bundle with logarithmic connection on $\CC$
   which we view as the other branch of $C$
   with coordinate $p$: 
    \[
   \Ec_{\hat H} \,\,=\,\,\biggl( E\otimes\Oc_{\CC}, \,\nabla = d+ (uv)) {d p\over p}\biggr). 
   \]  
   
   This means that the data $(\Ec_H, \Ec_{\hat H}, u,v)$ form an object of $\Mc_{\on{dR}}(C,\emptyset)$.
   So we get a functor $\Ien\to \Mc_{\on{dR}}(C,\emptyset)$. The fact that it is an equivalence, is verified in
   a standard way. \qed
   
   \vskip .2cm
   
   Theorem \ref {thm:micro-dr} is now obtained by gluing together the descriptions given by Proposition
   \ref{prop:can} over $X_{\on{sm}}$ and by Proposition \ref{prop:I} near the nodes of $X$. \qed

 \vfill\eject

  \section{Twisted microlocal sheaves}\label{sec:twisted-micro}
  
  \sparagraph{\bf Motivation: twisted $\Dc$-modules and sheaves.}
Let $X$ be a smooth algebraic variety over $\CC$. 
We recall \cite{BB} that to each class
 $t\in H^1_{\on{Zar}}
 \bigl(X, \bigl\{\Omega^1_X \buildrel d\over \to \Omega^{2, \on{cl}}_X\bigr\}\bigr)$
there corresponds a sheaf of rings of {\em twisted differential operators} on $X$
which we denote $\Dc_X^t$. 

Recall further that the first Chern class can be understood as a homomorphism
\[
c_1: \on{Pic}(X) \lra H^1_{\on{Zar}}
 \bigl(X, \bigl\{\Omega^1_X \buildrel d\over \to \Omega^{2, \on{cl}}_X\bigr\}\bigr).
\]
If $\Lc$ is a line bundle on $X$,  then we have an explicit model:
\[
\Dc_x^{c_1(\Lc)}  \,\,=\,\,\on{Diff}(\Lc, \Lc)
\]
is the sheaf formed by differential operators from sections of $\Lc$ to sections 
of $\Lc$. For a compact $X$, the  image of $c_1$ is typically an integer lattice 
in a complex vector space and
 the sheaves $\Dc_X^t$ can be seen as interpolating between the $\on{Diff}(\Lc, \Lc)$
 for different $\Lc$. We recall a particular explicit instance of this interpolation.
 
 Given a line bundle $\Lc$ on $X$, we denote by $\Lc^\circ$ the total space of $\Lc$
 minus the zero section, so $p:\Lc^\circ\to X$ is a $\CC^*$-torsor over $X$. 
 We denote by $\theta$ the Euler vector field ``$x \partial/\partial x$" on $\Lc^\circ$,
 i.e., the infinitesimal generator of the $\CC^*$-action. Thus $\theta$ is a global
 section of $\Dc_{\Lc^\circ}$. 
 
 \begin{prop}\label{prop:lambda}
 Let $\lambda\in\CC$. Then 
 \[
 \Dc_X^{\lambda c_1(\Lc)} \,\,\simeq \,\, p_*\biggl(
 \Dc_{\Lc^\circ} \biggl/  \Dc_{\Lc^\circ}(\theta-\lambda) \Dc_{\Lc^\circ}\biggr).
 \qed
 \]
 \end{prop}
 
 We now discuss the consequences of Proposition \ref{prop:lambda}
 for the Riemann-Hilbert correspondence for twisted $\Dc$-modules. 
 
 \vskip .2cm
 
On the $\Dc$-module side, the concepts of holonomic and regular $\Dc_X^t$-modules are defined in
the same way as in the untwisted case. We denote by $\Dc_X^{t}-\Mod_{\on{h.r.}}$ the category of 
holonomic regular $\Dc_X^t$-modules,
and by  $D^b_{\on{h.r.}}(\Dc_X^{t}-\Mod)$
the derived category formed by complexes with holonomic regular cohomology modules. 

\vskip .2cm

On the sheaf side,  choose $q\in\k^*$. Let $\Lc$ be a line bundle on $X$.
We denote by $\Sh^{\Lc, q}(X)$ the category of sheaves on $\Lc^\circ$ whose restriction on each fiber of $p$
is a local system with scalar monodromy $q\cdot\on{Id}$. Let $D^b(X)^{\Lc,q}$ be the bounded derived category
of $\Sh^{\Lc, q}(X)$. We denote by $\Dbc(X)^{\Lc, q} \subset D^b(X)^{\Lc, q}$ the full subcategory formed by complexes
with  $\CC$-constructible cohomology sheaves, and $\Perv^{\Lc, q}(X)\subset \Dbc(X)^{\Lc, q}$ the full subcategory of perverse sheaves.

 Proposition \ref{prop:lambda} implies the following. 
 
 \begin{cor} Take the base field $\k=\CC$.
  Let $\Lc$ be a line bundle on $X$ and $\lambda\in\CC$. 
 We have an  anti-equivalence of (pre-)triangulated categories and a compatible anti-equivalence of abelian categories
 \[
 D^b_{\on{h.r.}}(\Dc_X^{\lambda c_1(\Lc)}-\Mod) \to \Dbc(X)^{\Lc, e^{2\pi i \lambda}}, \quad 
  \Dc_X^{\lambda c_1(\Lc)}-\Mod_{\on{h.r.}}\to \Perv^{\Lc, e^{2\pi i\lambda}}(X).
 \quad
 \]
  
  \end{cor}

 \begin{rem} 
 For example, if $\lambda=n$ is an integer, then the monodromy comes out
 to be trivial, and we get that
 $\Dc_X^{\lambda c_1(\Lc)}-\Mod_{\on{h.r.}}$ is anti-equivalent
  to $\Perv_X$. This can also be seen directly, as $D_X^{n c_1(\Lc)}=
  \on{Diff}(\Lc^{\otimes n}, \Lc^{\otimes n})$ and so we have the ``solution
  functor" associating to any module $\Mc$ over this sheaf of rings  the
  complex
  \[
  \on{Sol}(\Mc) \,\,=\,\, R\Homu_{\on{Diff}(\Lc^{\otimes n}, \Lc^{\otimes n})}
  (\Mc, \Lc^{\otimes n}). 
  \]
  This complex is perverse, and the functor $\on{Sol}$ establishes the desired
  anti-equivalence. 
  
  \end{rem}
  
  We will also consider the ``universal twist" situation by not requiring the monodromy to
  be a fixed scalar multiple of 1 and working instead with monodromic sheaves and
  complexes on $\Lc^\circ$.
  
  That is, we consider the derived category
  $D^b_{\on{mon}}(\Lc^\circ)$ defined as the full subcategory in $D^b\Sh(\Lc^\circ)$
  formed by   $\CC$-monodromic complexes. 
  Inside it, let
  $\Dbc(X)^\Lc$ be the full triangulated subcategory of $\CC$-constructible $\CC$-monodromic complexes
  and $\Perv(X)^\Lc$ the abelian subcategory of  perverse sheaves on $\Lc^\circ$ which are
  $\CC$-monodromic. 
  
  Note that the natural functor $D^b (X)^{\Lc, q}\to  \Dbc(X)^\Lc$
  is not fully faithful. In the $\Dc$-module picture this correponds to the fact that
  the derived  pullback functor on modules corresponding to the projection of sheaves of rings
  $\Dc_{\Lc^\circ}\to \Dc_{\Lc^\circ}/(\theta-\lambda)$   is not fully faithful.

  \vskip .2cm
  
  \sparagraph {\bf Twisted microlocal sheaves.} 
  We now  modify the above and apply it to the case when $X$ is a nodal
  curve.  
  
  So let  $X$ be a nodal curve over $\CC$ with the normalization map
  $\varpi: \wt X\to X$, as in \S \ref{sec:microsheaves}.
  We denote by $D\subset X$ the set of nodes, and by $\wt D\subset X$ its
  preimage under $\varpi$. For any node $x$ we choose a small analytic neighborhood
   $U=U_x = B'\cup_x B''$ of $x$. Here
   $B', B''$ are  two branches of $X$ near $x$ which we identify with their preimages
  $\wt B', \wt B''\subset\wt X$. 
  
  Let $\Lc$ be a line bundle on $X$.
  We denote by $\wt \Lc = \varpi^{*}(\Lc)$ its pullback to $\wt X$ and by $\wt p: \wt\Lc^\circ\to\wt X$
  the projection. 
   For each node $x$ we choose an 
   {\em almost-trivialization}    of $\Lc$  over $U_x$,
   by which we mean an identification of $\Lc|_{U_x}$ with the trivial line bundle 
 with fiber   $\Lc_x$ or, equivalently, an
 identification of $\GG_m$-torsors
  \be\label{eq:almost}
  \Lc^\circ|_{U_x} \lra U_x\times\Lc_x^\circ
  \ee 
(Note that the space of almost-trivializations
 is contractible.)
 The isomorphism  \eqref{eq:almost}
   gives rise to the {\em relative},  or (fiberwise with respect to the projection to $\Lc^\circ_x$) Fourier transforms
  which are quasi-inverse equivalences of triangulated categories
  \[
  \xymatrix{
   D^b(\wt{B}', x')^{\wt\Lc}  
   \ar@<.7ex>[r]^{\on{FT}'}&  \ar@<.7ex>[l]^{\on{FT}''}
    D^b(\wt{B}'', x'')^{\wt\Lc},
  }
  \quad\quad
    \xymatrix{
   D^b(\wt{B}', x')^{\wt \Lc, q}  
   \ar@<.7ex>[r]^{\hskip -.7cm \on{FT}'}&  \ar@<.7ex>[l]^{\hskip -.7cm \on{FT}''}
    D^b(\wt{B}'', x'')^{\wt \Lc, q}, \,\,  q\in\k^*. 
  }
  \]
  They induce similar equivalences of abelian categories of twisted perverse sheaves.

\begin{Defi}\label{def:tw-micro} Let $q\in\k^*$.
 
(a) 
An $\Lc$-twisted, resp. $(\Lc, q)$-twisted microlocal complex on $X$ is a datum $\Fc$
consisting of:
\begin{enumerate}
\item[(1)] 
An object $\wt\Fc^\circ$ of $\Dbc(\wt X)^{\wt\Lc}$, resp. of $\Dbc(\wt X)^{\wt\Lc, q}$

\item[(2)]
For each node $x\in D$ with the two branches $B', B''$ as above,
isomorphisms
\[ 
\FT'(\wt\Fc^\circ |_{\wt p^{-1}(\wt B')}) \lra \wt\Fc^\circ |_{\wt p^{-1}(\wt B'')}, \quad 
\FT''(\wt\Fc^\circ |_{\wt p^{-1}(\wt B'')}) \lra \wt\Fc^\circ |_{\wt p^{-1}(\wt B')},
\]
inverse to each other. 
 \end{enumerate}

\vskip .2cm

\noindent (b) An $\Lc$-twisted, resp.  $(\Lc, q)$-twisted microlocal sheaf is an 
$\Lc$-twisted, resp. $(\Lc, q)$-twisted
microlocal
complex such that $\wt \Fc^\circ$ is a perverse sheaf on
 $\wt\Lc^\circ$.

\end{Defi}

 As before, for any finite subset $A\subset X$ of smooth points we denote
 by  $D\Mc^{\Lc}(X,A)$, resp.   $D\Mc^{\Lc, q}(X,A)$ the pre-triangulated
 dg-category formed by  $\Lc$-twisted, resp. $(\Lc, q)$-twisted microlocal complexes $\Fc$ on $X$
  such that $\wt\Fc^\circ$ has locally constant cohomology outside
  of the preimage of $A$ in $\wt\Lc^\circ$. By $\Mc^\Lc(X,A)$, resp.  $\Mc^{\Lc, q}(X,A)$
  we denote the full (abelian) subcategory in   $D\Mc^{\Lc}(X,A)$, resp.  $D\Mc^{\Lc, q}(X,A)$
  formed by $q$-twisted microlocal sheaves.
  
  \sparagraph{Calabi-Yau properties.} 
  Theorem \ref{thm:CY} generalizes to the twisted case as follows. 
  
  \begin{thm}\label{thm:CY-twist}
  Assume $X$ is a compact nodal curve with a duality structure, and
  $(X_i)_{i\in I}$ be its irreducible components.  Let $\Lc$ be a line bundle
  on $X$ with an  almost-trivialization on a neighborhood of each node. Then:
  \begin{enumerate}
  \item[(a)] $D\Mc^\Lc(X,\emptyset)$ is a Calabi-Yau category of dimension 3.
  
  \item[(b)] For any $q\in\k^*$ 
    we have that $D\Mc^{\Lc, q}(X,\emptyset)$
  is a Calabi-Yau category of dimension 2.
  
   \end{enumerate} 
  \end{thm}
  
  \begin{ex}
  If $X$ is a smooth projective curve of genus $g$, then part (a) corresponds to 
  the Poincar\'e duality on the compact 3-manifold $\Lc^\circ/\RR^*_+$,
  the circle bundle on $X$ associated to $\Lc$. 
  \end{ex}

  \noindent{\sl Sketch of proof of Theorem \ref {thm:CY-twist}:}
 It is obtained by arguments similar to those for 
  Theorem \ref{thm:CY}. That is, for any two objects $\Fc, \Gc$  of the category $D\Mc^\Lc(X,\emptyset)$
  resp. $D\Mc^{\Lc,q}(X,\emptyset)$ we introduce a constructible complex $\Mc\Homu^\Lc(\Fc,\Gc)$
 resp. $\Mc\Homu^{\Lc,q}(\Fc,\Gc)$ whose complex of global sections over $X$ is identified with
 $R\Hom(\Fc,\Gc)$ in the corresponding category.
  The statement then follows from canonical identifications
 \[
 \begin{gathered}
 \DD \Mc\Homu^\Lc(\Fc,\Gc) \,\,\simeq \,\,\Mc\Homu^\Lc(\Gc,\Fc)[3],
 \\
 \DD \Mc\Homu^{\Lc,q}(\Fc,\Gc) \,\,\simeq \,\, \Mc\Homu^{\Lc, q}(\Gc,\Fc)[2]. 
 \end{gathered}
 \]
These identifications are obtained by comparing the bifunctor $\Mc\Homu^\Lc$ with the
bifuctor $\mu\Homu$ of \cite{KS} applied to constructible complexes on
manifolds of the form
$\Lc^\circ|_Z$, where $Z$ is a unibranched subset of $X$. 
\qed

\vfill\eject

  \vfill\eject
  
  \section{ Multiplicative preprojective algebras
  }\label{sec:prepro}

  \sparagraph {\bf The definitions.} 
  We recall the definition of multiplicative preprojective algebras, following
  \cite{CBS} \cite{yamakawa}.

  \begin{conv}\label{conv:alg-cat}
 There is a very close correspondence between:
  \begin{itemize}
  \item[(1)]  $\k$-linear categories $\Cc$ with finitely many objects.
  
  \item[(2)] Their {\em total algebras}
 \[
  \Lambda_\Cc   \,\,=\,\,\bigoplus_{x,y\in\Ob(\Cc)}\Hom_\Cc(x,y).  
  \]
  \end{itemize}
  For instance,  each object $x\in\Cc$ gives an idempotent $\1_x\in\Lambda_\Cc$, 
 left $\Lambda_\Cc$-modules are the same
   as covariant  functors 
   $\Cc\to\on\Vect$, and so on.
    For this reason we will not make a notational distinction
   between objects of type (1) and (2), thus, for example,
    speaking about objects of an algebra   $\Lambda$ and morphisms between them
   (meaning objects and morphisms of a category $\Cc$ such that $\Lambda=\Lambda_\Cc$). 
  
  \end{conv}

  Let $\Gamma$ be a  {\em quiver},
   i.e., finite oriented graph, with the set of vertices $I$ and the set of
  arrows $E$, so we have the source and target maps $s,t: E\to I$. 
We fix a total ordering $<$ on $E$.   
  
  \begin{Defi}
  Let $\ul q=(q_i)_{i\in I}\in(\CC^*)^I$. The {\em multiplicative preprojective algebra}
   $\Lambda^{\ul q}(\Gamma)$
   is  
  defined by generators and relations as follows:
  \begin{enumerate}
  \item[(0)] $\Ob(\Lambda^{\ul q}(\Gamma)) = I$. In particular, for each $i\in I$ we have the
  identity morphism $\1_i: i\to i$. 
  
  \item[(1)] For each arrow $h\in E$ there are two generating morphisms
  $a_h: s(h)\to t(h)$ and $b_h: t(h)\to s(h)$.   We impose the condition that 
  \[
  \1_{t(h)} + a_h b_h: t(h)\to t(h), \quad \1_{s(h)} + b_h a_h: s(h)\to s(h)
  \]
  are invertible, i.e., introduce their formal inverses. 
  
  \item[(2)] We  further impose
  the following relations: for each $i\in I$,
  \[
  \prod_{h\in E: t(h)=i} (\1_i + a_h b_h) \prod_{h\in E, s(h)=i} (\1_i + b_h a_h)^{-1} 
  \,\,=\,\, q_i \1_i,
  \]
  where the factors in each product are ordered using the chosen total order 
  $<$ on $E$. 
   \end{enumerate}

   \end{Defi}

   It was proven in \cite[Th. 1.4]{CBS}
   that up to an isomorphism, $\Lambda^{\ul q}(\Gamma)$ is independent on the
   choice of the  order $<$, as well as on the choice of orientation of edges of $\Gamma$. 
   
   \vskip .2cm
   
   \sparagraph {\bf Microlocal sheaves on rational curves.} 
   Let now  $X$ be a  compact nodal curve over $\CC$ with the set of components $X_i, i\in I$.
   We then have the {\em intersection graph} $\Gamma_X$ of $X$.
   By definition, this is an un-oriented graph with the set of vertices
   $I$ and as many edges from $i$ to $j$ as there are intersection points of $X_i$ and
   $X_j$. In particular, for $i=j$ we put as many loops as there are self-intersection
   points of $X_i$. We now choose an orientation of $\Gamma_X$  and an ordering
   of the arrows in an arbitrary way, thus making it into a quiver, 
   so that the above constructions apply to $\Gamma_X$.
   Note that an  orientation of $\Gamma_X$ is the same as an
   orientation of  $X$ in the sense of  \S \ref {sec:dr}A. 
   
   Let $\Lc$ be a line bundle on $X$. We keep the notation of \S \ref{sec:twisted-micro}.
   Let $d_i = \deg(\varpi^{*}_i\Lc)\in\ZZ$. For $q\in\k^*$ we denote
   $q^{\deg(\Lc)}= (q^{d_i})_{i\in I}$. 
   
   \begin{thm}\label{thm:prepro} Assume
    that all the components $X_i$ are rational, i.e., the normalizations
    $\wt X_i$ are isomorphic to $\PP^1$. Then the category $\Mc^{\Lc, q}(X,\emptyset)$
    is equivalent to the category of finite-dimensional modules over
   $\Lambda^{q^{\deg(\Lc)}} (\Gamma_X)$. 
   
   \end{thm}

   \sparagraph {\bf Perverse sheaves on a disk:
   the $(\Phi, \Psi)$-description.}      The proof of Theorem
   \ref{thm:prepro} is based on a conceptual interpretation of the factors
    $\1_i + a_h b_h$
   and $(\1_i+ b_h a_h)^{-1}$ entering the defining relations of
   $\Lambda^{\ul q}(\Gamma)$. We observe that such expressions describe
    the  {\em monodromies of perverse sheaves on a disk}. 
    
    More precisely, let $B$ be an open  disk in the complex plane containing a point $y$.
    Let  $\overline B$ be an ``abstract"  closed disk containing $B$
    as its interior. 
    We denote $\Perv(B,y)$ the category of perverse sheaves on $B$ smooth
    everywhere except possibly $y$. Note that  for any $\Fc\in\Perv(B,y)$, 
    the restriction of $\Fc$ to $B-\{y\}$ is 
    a local system in degree 0  and so extends, by direct image,
    to a local system on $\overline B-\{y\}$.  So we can think of $\Fc$
    as a complex of sheaves on $\overline B$,  whose restriction to
    $\overline B-\{y\}$ is quasi-isomorphic to a local system in degree 0.
    In particular, for each $z\in \overline{B}-\{y\}$ we have a single vector
    space $\Fc_z$, the stalk of $\Fc$ at $z$.

     We have the following classical result
   \cite{beil-gluing} \cite{galligo-GM}. 
    \begin{prop}\label{prop:GGM}
   (a) Let $\Jen$
    be the category of diagrams of
    finite-dimensional $\k$-vector spaces
   \[
      \xymatrix{
 \Phi
   \ar@<.5ex>[r]^a&  \ar@<.5ex>[l]^b
  \Psi
  }
    \]
    such that the operator $T_\Psi=\1_\Psi+ ab$ is invertible. 
    For such a diagram the operator $T_\Phi = \1_\Phi+ba$ is invertible as well. 
     The category $\Perv(B,y)$ is equivalent to $\Jen$. 
    
    \vskip .2cm
    
    (b) Explicitly, an equivalence in (a) is obtained by choosing a boundary point
    $z\in\partial B$ and joining it with a simple arc $K$ with $y$. After such choices
    the vector spaces corresponding to $\Fc\in\Perv(B,y)$ are found as
    \[
    \Psi=\Psi(\Fc) = \Fc_z = \HH^0(K-\{y\}, \Fc), \quad \Phi = \Phi(\Fc) = \HH^1_K(B,\Fc).
    \]
    The operator $T_\Psi$ is the anti-clockwise monodromy of the local system
    $\Fc|_{B-\{y\}}$ around $y$. \qed
    \end{prop}
    
  The space $\Psi(\Fc)$ and $\Phi(\Fc)$ are referred to as the spaces
  of {\em nearby} and {\em vanishing cycles} of $\Fc$ at $y$ (with respect to the choice
  of an arc $K$).  
    
\vskip .2cm

\sparagraph {\bf Fourier transform in the $(\Phi, \Psi)$-description.}
Let $L$ be a 1-dimensional $\CC$-vector space, $L^*=\Hom_\CC(L,\CC)$
be the dual space, with the canonical pairing
\[
(z,w)\mapsto \langle z,w\rangle: L\times L^*\lra\CC.
\] 
Let $K$ be a half-line in $L$ originating at $0$, and
\[
K^* \,\,=\,\,\bigl\{ w\in L^*: \,\, \langle z,w\rangle \in\RR_{\geq 0}, \,\,
\forall z\in K\bigl\}
\]
be the dual half-line in $L^*$. 
We can consider $K$  as a simple arc  in $L$  
joining $0$ with the infinity of $L$, and similarly with $K^*$. Therefore
the choices of $K$ and $K^*$ give
identifications of the categories $\Perv(L,0)$ and $\Perv(L^*, 0)$
with the categories of diagrams as in Proposition \ref{prop:GGM}.

\begin{prop}\label{prop:cluster}
Under the identifications of Proposition \ref{prop:GGM},
the Fourier-Sato transform
\[
\FT: \Perv(L,0) \lra \Perv(L^*, 0)
\]
corresponds to the functor $\FT_\Jen$ which takes
\[
\bigl\{ \xymatrix{
 \Phi
   \ar@<.5ex>[r]^a&  \ar@<.5ex>[l]^b
  \Psi
  } \bigr\}  \,\, \buildrel \FT_\Jen\over\longmapsto \,\,
\bigl\{  \xymatrix{
 \Psi
   \ar@<.5ex>[r]^{a'}&  \ar@<.5ex>[l]^{b'}
 \Phi    
 } \bigr\},
\]
where $(a',b')$ are related to $(a,b)$ by the ``cluster transformation"
\[
 \begin{cases}
  a'=-b, \\
  b' = a(\1_\Phi + ba)^{-1}.
  \end{cases} \qed
\]
\end{prop}

\begin{cor}\label{cor:inverse}
In the situation of Proposition \ref{prop:cluster} we have
\[
\1+a'b'\,\, =\,\,  (\1+ba)^{-   
1}. 
\]
\end{cor}  

Note 
  that  this corollary prevents us from having
a naive  statement of the kind  ``Fourier transform interchanges
$\Phi$ with $\Psi$ and $a$ with $b$". 

\vskip .2cm

\noindent {\sl Proof of Proposition \ref{prop:cluster}:} 
We first establish  the identifications
\be\label{eq:phi-ft-psi}
\Psi(\FT(\Fc)) \,\,\simeq \,\,\Phi(\Fc). 
\ee
Let $K^\dagger\subset L$ be the half-plane
formed by $z$ such that $\Re \langle z,w\rangle \geq 0$ for each $w\in K^*$.
From the definition of $\FT$,  see \cite{KS}, \S 3.7 and the fact that $\Fc$ is $\CC^*$-monodromic,
we see  that $\Psi(\FT(\Fc))$, i.e., the stalk of $\FT(\Fc)$ at a generic point of the ray $K^*$,
is equal to the vector space $\HH^1_{K^\dagger}(L,\Fc)$.  But $K^\dagger$ contains $K$ and can
be contracted to it without changing the cohomology with support for any $\Fc\in\Perv(L,0)$.
This means that $\Psi(\FT(\Fc)) \simeq \HH^1_{K}(L,\Fc) = \Psi(\Fc)$.

 \vskip .2cm

 Next, we prove the  Corollary \ref{cor:inverse}. 
Note that rotating $K$ in $L$ anti-clockwise
results in rotating $K^*$ in $L^*$ clockwise.   So the monodromy on
$\Phi(\Fc)$ obtained by rotating $K$ in the canonical way given by
the complex structure (i.e., anti-clockwise),   is the inverse of the
monodromy on $\Psi(\FT(\Fc))= \Phi(\Fc)$ obtained by rotating
$K^*$ in the same canonical way (i.e., also anti-clockwise). 
This establishes the  corollary. 

\vskip .2cm

We now prove  Proposition \ref{prop:cluster} in full generality
by using the approach of  \cite{beil-gluing}. 
 We identify $\Perv(L,0)$ with $\Jen$ throughout. Note that $m=(T_\Phi, T_\Psi)$
defines an automorphism of the identity functor of $\Jen=\Perv(L,0)$ called the
{\em monodromy operator}. Further, $\Perv(L,0)$ splits into
a direct sum of abelian categories
\[
\Perv(L,0) \,\,=\,\, \Perv(L,0)_u\oplus \Perv(L,0)_{nu}.
\]
Here $m$ acts unipotently on every object $\F\in \Perv(L,0)_u$
(equivalently, on $\Phi(\F)$, $\Psi(\F)$ for $\F\in \Perv(L,0)_u$), while $\1-m$ is invertible
on every object $\F\in \Perv(L,0)_{nu}$ .

We construct the isomorphism claimed in Proposition 
\ref{prop:cluster} separately for $\F\in \Perv(L,0)_u$ and $\F\in \Perv(L,0)_{nu}$.

\vskip .2cm

Assume first that $\F\in \Perv(L,0)_{nu}$. Notice that for $\F\in \Perv(L,0)_{nu}$
the maps $a:\Phi(\F)\to \Psi(\F)$ and $b:\Psi(\F)\to \Phi(\F)$ are invertible.
 This means that either
of the two functors $\F\mapsto (\Psi(\F), T_\Psi)$, $\F\mapsto (\Phi(\F), T_\Phi)$ is an equivalence
between $\Perv(L,0)_{nu}$ and the category of vector spaces with an automorphism which does 
not have eigenvalue  one. Thus in this case it suffices to construct a functorial isomorphism 
$\Phi(\F)\cong \Psi(\FT (\F))$ sending the automorphism $T_\Phi$ to $T_\Psi^{-1}$.
This reduces to Corollary \ref{cor:inverse}. 

\vskip .2cm

We now consider $\F\in  \Perv(L,0)_u$. Notice that the category $\Perv(L,0)_u$
has, up to isomorphism, two irreducible objects,  $\LL_0 = \ul\k_0[-1]$
and $\LL_1=\ul\k_L$  (the sky-scraper at zero and  the constant sheaf). 
Let  $\Pi_0$,
$\Pi_1$ be  projective covers of $\LL_0, \LL_1$, which are projective
objects in the category of pro-objects
\[
\on{Pro} \bigl( \Perv(L,0)_u\bigr) \,\,\subset \,\, \on{Fun}\bigl( \Perv(L,0)_u,
\Vect\bigr)^{\on{op}}. 
\]
 They are defined uniquely up to an isomorphism.
Moreover, any exact functor from $ \Perv(L,0)_u$ to vector spaces
sending $\LL_1$ (resp. $\LL_0$) to zero and $\LL_0$ (resp. $\LL_1$) 
to a one dimensional space is isomorphic, in the sense of viewing pro-objects as functors
above, 
to $\Pi_0$ (resp. $\Pi_1$). This means that there exist isomorphisms 
of functors $\Perv(L,0)\to\Vect$
\[
\Hom(\Pi_0, - ) \,\,\cong \Phi,\quad \Hom(\Pi_0, - )\,\,\cong \,\, \Psi. 
\]
  We fix such isomorphisms.

Proposition \ref{prop:GGM}
 implies that $\on{End}(\Pi_0)\simeq  \kb [[(m-1)]]\simeq  \on{End}(\Pi_1)$ while each of the
 spaces $\Hom(\Pi_0, \Pi_1)$, $\Hom(\Pi_1,\Pi_0)$ is a free rank one module over 
 $\kb [[(m-1)]]$ generated respectively  by elements $a$, $b$. 
 
 Since $\FT$ interchanges $\LL_0$ and $\LL_1$, we have
 \be\label{eq:phi-psi-pi}
 \FT(\Pi_0)\,\,\simeq \,\,  \Pi_1, \quad 
 \FT(\Pi_1)\,\,\simeq \,\, \Pi_0.
 \ee
 Furthermore, the isomorphism   \eqref{eq:phi-ft-psi}
 sending 
  $m_\F$ to $m_{\FT(\F)}^{-1}$ shows that for some (hence for any) choice of the isomorphisms
 $\FT(\Pi_0)\cong \Pi_1$ the automorphim $\FT(m)$ of the left hand side
 corresponds to the automorphism $m^{-1}$ of the right hand side. It follows that 
 an isomorphism $\FT(\Pi_1)\cong \Pi_0$ also sends $\FT(m)$ to $m^{-1}$. We can choose
 the isomorphisms  \eqref{eq:phi-psi-pi} in such a way that the map $\FT(a)$ becomes compatible with $-b$. 
 This is clear since both elements generate the corresponding free rank one modules over
 $\kb [[(m-1)]]$. Then we see that $\FT(b)$ corresponds to $a(\1_{\Pi_0}+ba)^{-1}$, this
 implies the statement. \qed

\begin{rem} In the last paragraph of the proof we made a choice of  isomorphisms
 \eqref{eq:phi-psi-pi} satisfying certain requirements. 
 We have earlier  constructed  an  isomorpism of functors  \eqref{eq:phi-ft-psi}. 
 Combining
 it with the canonical isomorphism $\FT^2(\F)=(-1)^*(\F)$ we can (upon making a binary choice
 of a homotopy class of a path connecting the ray $K$ to the ray $-K$) produce a canonical
 isomorphism $\Psi(\F)\simeq \Phi(\FT(\F))$.  
  These two isomorphism of functors yield
  isomorphisms of representing objects. 
 We do not claim however  that  these isomorphisms  
  satisfy our requirements. They provide another (isomorphic but different) functor on the category
 of linear algebra data of Proposition \ref{prop:GGM}; it may be interesting to work it out
 explicitly. 
 \end{rem}
 
 \begin{rem}
 In  the case $\k=\CC$
one can deduce the proposition
 from the infinitesimal description $\Perv(\CC,0)\simeq \Ien$
(Proposition \ref{prop:I}), where the Fourier transform functor
$\FT_\Ien: \Ien\to \Ien$ is given by Corollary \ref{cor:FT-I}:
\be\label{eq:FT-I}
\bigl\{ \xymatrix{ E
   \ar@<.5ex>[r]^{u}&  \ar@<.5ex>[l]^{v}
 F}\bigr\} \,\,\longmapsto \,\, \bigl\{ \xymatrix{ E'=F
   \ar@<.5ex>[r]^{u'}&  \ar@<.5ex>[l]^{v'}
 F'=E}\bigr\} , \quad 
 u'=v, v'=-u. 
\ee
Since both $\Ien$ and $\Jen$ describe $\Perv(\CC,0)$,
we get an  identification $\Ien\to \Jen$ which was given explicitly  by Malgrange \cite[(II.3.2)]{Ma} as follows:
\be\label{eq:Mal-corr}
\begin{gathered}
\bigl\{ \xymatrix{ E
   \ar@<.5ex>[r]^{u}&  \ar@<.5ex>[l]^{v}
 F}\bigr\} \,\,\longmapsto \,\, \bigl\{ \xymatrix{ \Phi=
   \ar@<.5ex>[r]^{a}&  \ar@<.5ex>[l]^{b}
 \Psi=F}\bigr\} 
 \\
 \begin {cases}
 a=u,\\
 b= \varphi(vu)\cdot v, \quad \varphi(z) = (e^{2\pi i z} -1)/z. 
 \end{cases}
 \end{gathered}
\ee
By inverting   \eqref{eq:Mal-corr} (i.e., finding $u$ and $v$ through $a$ and $b$),
and then applying  \eqref{eq:Mal-corr}  to $u', v'$ given by 
 \eqref{eq:FT-I}, we get an  object
 $\bigl\{ \xymatrix{ \Psi
   \ar@<.5ex>[r]^{a^\dagger}&  \ar@<.5ex>[l]^{b^\dagger}
 \Phi}\bigr\} $
 which turns out to be isomorphic to $\bigl\{  \xymatrix{
 \Psi
   \ar@<.5ex>[r]^{a'}&  \ar@<.5ex>[l]^{b'}
 \Phi    
 } \bigr\}$ by conjugation with an explicit invertible function of $\1_\Phi+ba)$.
 \end{rem}

\vskip .2cm

\sparagraph {\bf Proof of Theorem \ref{thm:prepro}}.  We start with an
almost obvious model
case of one projective line $Y\simeq \PP^1$. Suppose we are given
 a point $z\in Y$
which will serve as an ``origin" and a further set of 
$N$ points $A=\{y_1, \cdots, y_N\}$ which we position on the boundary of
a closed disk $B$ containing $z$, in the clockwise order.
   Choose a system of simple arcs $K_\nu$ joining $z$ with $y_\nu$
   and not intersecting outside of $z$. Let $\Lc $ be a line bundle
   of degree $d$ on $\PP^1$ and let $q\in\k^*$.
   
   \begin{lem}\label{lem:twisted}
   The category $\Perv^{(\Lc, q)}(Y,A)$ is equivalent to the category
   of diagrams consisting of vector spaces $\Psi, \Phi_1, \cdots, \Phi_N$
   and linear maps
   \[
   \bigl\{ \xymatrix{
 \Phi_\nu
   \ar@<.5ex>[r]^{a_\nu}&  \ar@<.5ex>[l]^{b_\nu}
  \Psi
  } \bigr\}, \quad \nu=1,\cdots, N, 
   \]
   such that
   each $\1_\Psi + a_\nu b_\nu$ is invertible and
   \[
   \prod_{\nu=1}^N (\1_\Psi+a_\nu b_\nu) \,\,=\,\, q^d \1_\Psi. 
   \]
   \end{lem}
   
   \noindent {\sl Proof:} We first consider the untwisted case: $q=1$ or, equivalently,
   no $\Lc$. In this case the statement follows at once from Proposition \ref{prop:GGM}. 
   Indeed, choose thin neighborhoods $U_\nu$ of $K_\nu$
   (thus containing $z$ and $y_\nu$ which are topologically disks and let $U=\bigcup U_\nu$.
   We can assume that $Y$ is, topologically, a disk as well.  An object
   $\Fc\in\Perv(Y,A)$ can be seen
   as consisting of  perverse sheaves $\Fc_\nu$ on $U_\nu$
  which are glued together into
 a global perverse sheaf on $Y$. Each $\Fc_\nu$ is described by
 a diagram
 $
 \bigl\{ \xymatrix{
 \Phi_\nu
   \ar@<.5ex>[r]^{a_\nu}&  \ar@<.5ex>[l]^{b_\nu}
  \Psi_\nu
  } \bigr\},$
  To glue the $\Fc_\nu$ together, we need, first, to identify all the $\Psi_\nu$
  with each other, i.e., with a single vector space $\Psi$.
   This will give a perverse sheaf  $\Fc_U$ on  $U$.
  In order for $\Fc_U$ to extend  to a perverse sheaf on $Y=\CC\PP^1$,
  it is necessary and sufficient that the monodromy of $\Fc_U$ along the
  boundary  $\partial U$ of $U$ be trivial, in which case the extension is unique up to
  a unique isomorphism. 
  
  To identify this condition explicitly,  let
   $\gamma_\nu$ be a loop in $Y$ beginning at $z$, going towards $y_\nu$ along $K_\nu$,
 then circling around $y_\nu$ anti-clockwise and returning back to $z$
 along the same path. 
 Then $\partial U$ can be represented, up to homotopy, by
  the composite loop 
 $\gamma = \gamma_1 \gamma_2\cdots\gamma_N$  
   and the monodromy of $\Fc_\nu$
 around $\gamma_\nu$ is $1+a_\nu b_\nu$.

 In the twisted case,  choose a trivialization of $\Lc$ over $U$, so that we have the
 projections
 \[
 U\buildrel \alpha\over\lla \Lc^\circ|_U \buildrel\beta\over\lra \CC^*
 \]
 Let $\wt z$ be the vector in the fiber of $\Lc$ over $z$ such that $\beta(\wt z)=1$. 
 Let 
 \[
 \wt\gamma = \gamma\times\{ 1\} \,\,\subset  U\times \CC^* \simeq \Lc^\circ|_U 
 \]
  be the lift of $\gamma$ with respect to the trivialization. Since $\gamma$
  does not meet $A$, we can regard $\wt\gamma$ as
  a loop in  $\Lc^\circ|_{Y-A}$, beginning and ending at $\wt z$.
 
 Note that the line bundle $\Lc$ is trivial over $Y-A$ as well, and so
 \[
 \pi_1\bigl( \Lc^\circ|_{Y-A}, \wt z\bigr) \,\,=\,\, \ZZ\cdot \zeta,
 \]
 where $\zeta$ is the counterclockwise loop in the fiber $\Lc^\circ|_z$. 
 Under this identification, the element represented by $\wt\gamma$,
 is equal to $d\cdot\zeta$. 
 
 Now,  using our trivialization,
 we have an equivalence
 \[
M:  \Perv(U,A) \lra \Perv^{(\Lc, q)}(U,A), \quad \Fc \mapsto \alpha^*\Fc \otimes_\k \beta^*\Ec_q,
 \]
 where $\Ec_q$ is the 1-dimensional local system on $\CC^*$ with monodromy $q$.
 An object $\Fc$  of $\Perv(U,A)$ is described by a diagram 
   of
 \[
  \bigl\{ \xymatrix{
 \Phi_\nu
   \ar@<.5ex>[r]^{a_\nu}&  \ar@<.5ex>[l]^{b_\nu}
  \Psi
  } \bigr\}, \quad \nu=1, \cdots, N
  \]
  as before. 
 The possibility of extending $M(\Fc)$ 
   from $\Lc^\circ|_U$ to the whole of
 $\Lc^\circ$ is equivalent to the  monodromy around $\zeta \in \pi_1(\Lc|_{Y-A},\wt z)$ being
 equal to $q\cdot \1$. In view of the equality $\wt\gamma = d\cdot\zeta$, this
  gives precisely the condition of the lemma. \qed

 \vskip .3cm
 
 The proof of  Theorem \ref{thm:prepro} is now obtained by gluing
 together  the descriptions of Lemma
 \ref{lem:twisted}, using Proposition \ref{prop:cluster} and Corollary
 \ref{cor:inverse}. 
 
 More precisely,  we apply the lemma to each $(Y_i, A_i)$, $i\in I$, where the
 $Y_i=\wt X_i\buildrel\varpi_i\over\to X$, $i\in I$ are the components of the normalization $\wt X$
 of $X$, and $A_i = \wt D\cap \wt X_i$.  We recall that
   $\wt D\subset \wt X$ is the preimage of the set of nodes $D\subset X$. 
   We put $\Lc_i = \varpi_i^*\Lc$ so that $d_i=\deg(\Lc_i)$. 
  
  Choose an orientation of the intersection graph $\Gamma = \Gamma_X$, or, equivalently, an  ordering $(x', x'')$
  of the pair of preimages of each node $x\in D$. We will label these preimages by the arrows $h$ of $\Gamma$,
  i.e.,  denote them by
  \[
  x'_h \in A_{s(h)} \subset Y_{s(h)}, \quad x''_h \in A_{t(h)} \subset Y_{t(h)}, \quad  h\in E. 
  \]
  Thus $A_i$ consists of
  \[
  x'_h, \,\, s(h) = i \text{ and } \,\, x''_{h}, \,\, t(h)=i.
  \]
  We choose a base point $z_i$ in each $Y_i$ and position the elements of $A_i$ on the boundary
  of a disk around $z_i$, so that, in the clockwise order, we have first the $x''_h, t(h)=i$
  (according to the order $<$ on $E$) and then the $x'_h, s(h)=i$ (again according to $<$). 
  We join $z_i$ with the elements of $A_i$ simple arcs meeting only at $z_i$. 
  
 An object $\Fc_i\in\Perv^{(\Lc_i, q_i)}(Y_i, A_i)$ is then described by a diagram consisting of one space
 $\Psi_i$ together with spaces $\Phi_{x'_h}$, $s(h)=i$ and $\Phi_{x''_h}$, $t(h)=i$ together with the maps
 \[
   \bigl\{ \xymatrix{
 \Phi_{x'_h}
   \ar@<.5ex>[r]^{a'_{h}}&  \ar@<.5ex>[l]^{b'_{h}}
  \Psi_i
  } \bigr\}, \,\, s(h)=i, \quad
    \bigl\{ \xymatrix{
 \Phi_{x''_h}
   \ar@<.5ex>[r]^{a''_h}&  \ar@<.5ex>[l]^{b''_h}
  \Psi_i
  } \bigr\},\,\,\, t(h)=i
 \] 
 so that the condition of the lemma reads:
 \be\label{eq:perv-product}
 \prod_{t(h)=i} (\1+ a''_h b''_h)  \prod_{s(h)=i} (\1+ a'_h b'_h) \,\,\,=\,\,\, q^{d_i}\cdot \1. 
 \ee
In order to glue the $\Fc_i$ into one
  twisted microlocal sheaf on $X$, we need to specify an identification of Fourier transforms
  at each node $x$. This means that $\Psi_i$ (which is identified with the space of nearby cycles
  of $\Fc_i$ at each $x'_h$, $s(h)=i$ and each $x''_h, t(h)=i$) becomes identified with
  the space of vanishing cycles of
  $\Fc_{t(h)}$ at $x''_h$ for $s(h)=i$ and of $\Fc_{s(h)}$ at $x'_h$ for $t(h)=i$. 
   
  Therefore all the linear algebra data reduce to the vector spaces $V_i=\Psi_i$ and
  linear operators
  \[
  \begin{gathered}
  a_h: V_{s(h)} = \Psi_{s(h)} \simeq \Phi_{x''_h}(\Fc_{t(h)}) \buildrel a''_h\over\lra \Psi_{t(h)} = V_{t(h)}, 
  \\
  b_h: V_{t(h)} = \Psi_{t(h)} \buildrel b''_h \over\lra \Phi_{x''_h}(\Fc_{t(h)} \simeq \Psi_{s(h)} = V_{s(h)},
  \end{gathered}
  \]
  where $\simeq$ stands for the identifications given by the  Fourier transform. 
  This means that we do not use the simply primed $a'_h, b'_h$, expressing them
  through $a''_h, b''_h$ by Proposition \ref{prop:cluster}. 
  
 After this reduction,
  the conditions  \eqref{eq:perv-product}  coincide, in view of Corollary \ref{cor:inverse}, with the defining relations
  of the multiplicative preprojective algebra. \qed
  
  \vfill\eject
  
  \section{Preprojective algebras for general nodal curves}\label{sec:preprogen}
  
 Theorem \ref{thm:prepro} can be extended  to the case of arbitrary compact nodal curves by introducing an appropriate
 analog of preprojective algebras (PPA).  In this section we present this analog and discuss possible further
 generalizations to differential graded (dg-) case and their consequences for the symplectic structure of
 moduli spaces. 
 
 Throughout the paper we use the notation
 \[
 [\alpha, \beta] \,\,=\,\, \alpha\beta \alpha^{-1}\beta^{-1}
 \]
  to denote the group commutator. 
  
  \sparagraph{Higher genus PPA.}
  Let $X$ be a  compact nodal curve over $\CC$.  As before we denote by 
   $D$ the set of nodes of $X$, by $X_i, i\in I$
   the irreducible components of $X$ and by
    $\wt X_i\subset \wt X \buildrel\varpi\over\to X$ the normalizations
   of $X_i$ and $X$. Let $\Lc$ be a line bundle on $X$ and $\wt\Lc = \varpi^*\Lc$. We denote by:
   \[
   \begin{gathered}
   g_i = \text{ the genus of   } \wt X_i,  \quad 
    d_i = \deg(\wt \Lc|_{\wt X_i}), \quad \wt D_i = \varpi^{-1}(D)\cap \wt X_i.  
   \end{gathered} 
   \]
   We choose an orientation of $X$, i.e., a total order $x' < x''$ on each
   2-element set $\varpi^{-1}(x), x\in D$,  see \S \ref {sec:dr}A. 
   
   For each node $x\in D$ we denote by $s(x)\in I$
   the label of the irreducible component containing $x'$, and by $t(x)$
   the label of the component containing $x''$. We also choose a total
   order on the set $D$.

     \begin{Defi}\label{def:hg-ppa}
   Let $X, \Lc$ as above be given and $q\in\k^*$. The 
   {\em preprojective $(X,\Lc)$-algebra} $\Lambda^{\Lc, q} (X)$  is defined
   by generators and relations as follows:
   
    \begin{enumerate}
  \item[(0)]   Objects $i\in I$.
  
  \item[(1)] For each node $x\in D$,  two generating morphisms 
 $a_x: s(x)\to t(x)$ and $b_x: t(x)\to s(x)$.   We impose the condition that 
  \[
  \1_{t(h)} + a_h b_h: t(h)\to t(h), \quad \1_{s(h)} + b_h a_h: s(h)\to s(h)
  \]
  are invertible, i.e., introduce their formal inverses. 
  
  \item[(1')] For each $i\in I$ there are generating morphisms  
 \[
  \alpha^i_\nu, \beta^i_\nu, \,\,\, i=1, \cdots, g_i, 
   \]
  which are required to be invertible. 
  
   \item[(2)]  For each $i\in I$ we impose a relation
 
 \[
   \begin{gathered}
 \prod_{x\in D: t(h)=i} (\1_i + a_x b_x) \prod_{x\in D, s(x)=i} (\1_i + b_x a_x)^{-1}  
  \prod_{\nu=1}^{g_i} [\alpha^i_\nu, \beta^i_\nu]
 \,\,= \,\, q^{d_i} \1_i.
  \end{gathered} 
  \]
  Here the factors in the first two products are ordered using the
   chosen total order 
  $<$ on $D$. 
   \end{enumerate} 
      \end{Defi}
      
         \begin{exas}  
    (a)   If   all $X_i$  are rational, then
      $\Lambda^{\Lc, q}(C)$ reduces to the multiplicative preprojective
      algebra associated to the quiver $\Gamma_X$, 
and parameters $q^{d_i}$,       
      see 
      \S \ref{sec:prepro}.   
      
      (b) If $X$ is smooth irreducible of genus $g>0$, then the fundamental group
      $\pi_1(X)$ has a universal central extension $\wt\pi_1(X)$ given by generators and relations
      as follows
      \[
   \wt\pi_1(X) \,\, = \,\,   \biggl\langle \alpha_1, \cdots, \alpha_g, \beta_1, \cdots, \beta_g, \q\, \biggl| 
    \,\,\, \prod_{\nu=1}^g [\alpha_\nu \beta_\nu] = \q, \,\, [\alpha_i,\q] = [\beta_i, \q]=1 \biggr\rangle. 
      \]
      In this case $\Lambda^{\Lc, q}(X)$ is
      a quotient of the group algebra of  $\wt\pi_1(X)$ by the relation $\q=q^d$. 
          \end{exas}

            \begin{thm}\label{thm:higher-pre}
  The abelian category   $\Mc^{\Lc,q}(X,\emptyset)$
      is equivalent to  the category of finite-dimensional modules over
      $\Lambda^{\Lc,q}(X)$.    \end{thm}
      
      \sparagraph{Proof of Theorem \ref{thm:higher-pre}.}
      It is similar to that of Theorem \ref{thm:prepro}. We first consider
      the following model case.
      
      Let $Y$ be a smooth, compact, irreducible curve of genus $g$  
      together with finite subset $A=\{y_1, \cdots, y_N\} 
      \subset Y$. Let $\Lc$ be a line bundle 
      over $Y$ of degree $d$. 
      Define a $\k$-algebra $\Lambda^{\Lc,q}(Y,A)$ by 
      generators and relations as follows;
      \begin{enumerate}
      \item[(0)] Objects $\psi$, $\phi_1, \cdots, \phi_N$.
      
      \item[(1)] Generating morphisms  
      \[
      \begin{gathered}
      a_\lambda: \phi_\lambda\to \psi, \quad b_\lambda: \psi\to\phi_\lambda,\quad 
      \lambda = 1,\cdots, N; 
      \\
      \alpha_\nu, \beta_\nu: \psi\to\psi, \,\,\, \nu=1, \cdots, g.
     \end{gathered}
      \]
      We  require that 
      \[
      \1_\psi + a_\lambda b_\lambda, \,\,\, 
       \1_{\phi_\lambda}+b_\lambda a_\lambda, \,\,\, a_\nu, b_\nu, h_\mu
       \]
        be invertible, i.e., introduce
       their formal inverses.
       
        \item[(2)]  One relation       \[
      \prod_{\lambda=1}^N (\1_\psi + a_\lambda b_\lambda) 
       \prod_{\nu=1}^g [\alpha_\nu, \beta_\nu]    \,\,=\,\, 
        q^d \1_\psi. 
       \]
      \end{enumerate}

     \begin{lem}\label{lem:higher-perv}
     
     The abelian category $\Perv^{\Lc}(Y, A)$ is equivalent to 
     the category of finite-dimensional $\Lambda^\Lc(Y,A)$-modules. 
     \end{lem}
     
    \noindent {\sl Proof:}   Completely analogous to that of Lemma \ref{lem:twisted}. We choose a base
    point $p\in Y-A$, realize $\alpha_i$ and $\beta_i$ as the standard A- and B-loops based at $p$
    and choose simple arcs $K_\lambda$ jointing $p$ with $y_\lambda$ so that they do not intersect
    except at $p$ and follow each other in the clockwise order. Conjugating with  $K_\lambda$  a small
    loop around $y_\lambda$, we get a loop $h_\lambda$ based at $p$, and we can choose
    the $K_\lambda$ to follow the system of  $\alpha_i,  \beta_i$ in the clockwise order so that   in $\pi_1(Y-A, p)$  
    we have the relation
   \[
   \prod_{\lambda=1}^N h_\lambda \prod_{\nu=1}^g [\alpha_\nu, \beta_\nu] \,\,=\,\, 1. 
   \]
    Let $D$ be a disk containing all the paths $K_\lambda$, so $\Lc$ is trivial over $D$.
    The lemma is obtained by gluing
    the category of  perverse sheaves on $D$   and that of (twisted) local systems on $X-D$. 
    \qed
   
Theorem \ref{thm:higher-pre} is now obtained by gluing the descriptions of Lemma \ref{lem:higher-perv} for
$(Y,A)=(\wt X_i, \wt D_i)$  for various $i$. \qed

\sparagraph{Remarks on derived PPA.} The algebra  $\Lambda^{\Lc,q}(X)$ has a derived analog.
This is a dg-algebra $L\Lambda^{\Lc,q}(X)$ with the same generators $a_x, b_x, \alpha^i_\nu, \beta^i_\nu$
as $\Lambda^{\Lc,q}(X)$ (considered in degree $0$)
 with the same conditions of invertibility but instead of imposing relations 
in Definition \ref {def:hg-ppa}, we introduce new free generators of degree $-1$ whose differentials
are put to be the differences between the LHS and RHS of the relations. The symbol $L$ is used to
signify the left derived functor. Thus $\Lambda^{\Lc,q}(X)$ is the $0$th cohomology algebra of 
  $L\Lambda^{\Lc,q}(X)$. 
  
  \vskip .2cm
  
  It seems very likely that the triangulated category $D\Mc^{\Lc,q}(X)$ can be identified with the
  derived category formed by finite-dimensional dg-modules over $L\Lambda^{\Lc,q}(X)$
  (with quasi-isomorphisms inverted). In view of Theorem \ref{thm:CY}  we can then  expect
  that  $D\Mc^{\Lc,q}(X)$ is a Calabi-Yau dg-algebra of dimension 2. In other word, we expect that,
  denoting $L = L\Lambda^{\Lc,q}(X)$,  there is  a quasi-isomorphism of $L$-bimodules 
  \be\label{eq:hoch} 
 \gamma: L\to  L^! :=  R\Hom_{L\otimes L^{\on{op}}}(L, L\otimes L^{\on{op}})[2], \quad
 \text{such that} \quad  \gamma = \gamma^![2],
  \ee
    see \cite{Gi}, Def. 3.2.3. 
  
   In general,   $L\Lambda^{\Lc,q}(X)$ is not quasi-isomorphic to  $\Lambda^{\Lc,q}(X)$,
  which explains the following example.

   \begin{ex}
      Let $X$ be the union of two projective lines meeting transversely,
      let $\Lc$ be trivial and $q=1$. Then
       $D\Mc(X,\emptyset)$ is a Calabi-Yau category of
     dimension 2,  while
      $\Mc(X,\emptyset)$
      has infinite cohomological dimension.
      Indeed, $\Mc(X,\emptyset)$
      is identified with the category of modules over the 
     multiplicative preprojective algebra
     corresponding to the quiver $A_2$; 
     this algebra
      has two objects $1,2$ generating morphisms $a:1\to 2$ and $b:2\to 1$
     subject to the relations $ab=ba=0$.
      \end{ex}
      
     We can also define the {\em universal higher genus PPA} (derived as well as non-derived)
     by replacing  $q\in\k^*$ in the above by an indeterminate $\q$ and working over the Laurent
     polynomial ring $\k[\qpm]$. We denote the corresponding (dg-) algebras by
     $L\Lambda^\Lc (X)$ and $\Lambda^\Lc(X)$. 
     
      Because of the 1-dimensionality of $\k[\qpm]$, we expect that  $L\Lambda^\Lc (X)$,
      considered as a dg-algebra over $\k$, 
      is   3-Calabi-Yau, rather than  2-Calabi-Yau. 
   
 \begin{ex}
 If $X$ is a smooth projective curve of genus $g>0$, then $\Lambda^\Lc(X)$
 is  the group algebra of the fundamental group of $\Lc^\circ$. Now, $\Lc^\circ$
 is homotopy equivalent to a circle bundle over $X$, which is a compact, apsherical,
 oriented 3-manifold.  By \cite{Gi}, Cor. 6.1.4 this implies that
 $\Lambda^\Lc(X)$ is a (non-dg) 3-Calabi-Yau algebra.   
 Further, in this case $L\Lambda^\Lc (X)$ is quasi-isomorphic to
  $\Lambda^\Lc(X)$ by \cite{Gi}, Thm. 5.3.1. 
 \end{ex} 
 
 \sparagraph{Remarks on moduli spaces.} Assume $\on{char}(\k)=0$. We would like to view the
 symplectic nature of (multiplicative) quiver varieties as yet another manifestation of the
 following general principle, which  also encompasses the approaches of \cite{goldman}
 and \cite{Mu} to local systems (resp. coherent sheaves) on topological (resp. K3 or abelian)
 surfaces.
 
 \begin{CYP}\label{CYP}
 If $\Cc$ is a Calabi-Yau category of dimension 2, then $\Men$, the ``moduli space"
 of objects in $\Cc$, has a canonical symplectic structure. After identifying the
 ``tangent space" to $\Men$ at the point corresponding to object $E$, with
 $\Ext^1_\Cc(E,E)$, the symplectic form is given by the {\em cohomological pairing}
 \[
 \Ext_\Cc^1(E,E) \otimes \Ext_\Cc^1(E,E) \buildrel\cup\over\lra \Ext^2_\Cc(E,E)
  \buildrel \on{tr}\over\lra
 \k,
 \]
 where $\on{tr}$ corresponds, via the Calabi-Yau isomorphism, to the embedding $\k\to\Hom_\Cc(E,E)$. 
 \end{CYP} 
 
 This principle, along with a generalization to $N$-Calabi-Yau categories for any $N$,
 was formulated in \cite{kontsevich-soibelman} \S 10 and made
  precise in a formal neighborhood context.
 A wider, more global,  interpretation would be as follows.
 
 \vskip .2cm
 
 \noindent {\bf ``Space":} understood in the sense of derived algebraic geometry
 \cite{lurie-DGA} \cite{TVe}, as a {\em derived stack}.  
Informally,  a derived stack $\Yen$ can be seen as a nonlinear (curved)  analog of
a cochain complex  of $\k$-vector spaces, in the same sense in which a manifold
can be seen as a curved analog of a single vector space. 
In particular, for a $\k$-point $y\in\Yen$ we have the {\em tangent dg-space}
$T^\bullet_y\Yen$, which is a cochain complex. The {\em amplitude} of $\Yen$ is
an integer interval $[a,b]$ such that 
$H^i T^\bullet_y\Yen=0$ for $i\notin [a,b]$ and all $y$. 
Given a morphism $f: Y\to Z$ of smooth affine algebraic varieties over $\k$ and a
$\k$-point $z\in Z$, we have the  {\em derived  preimage} $Rf^{-1}(z)$, which is 
a derived stack (scheme) of amplitude $[0,1]$, see \cite{CFK} for elementary treatment.

\vskip .2cm

 \noindent {\bf ``Moduli":} understood as the derived  stack $\Men_\Cc$
  of moduli of objects in
 a dg-category $\Cc$  defined in \cite{TV}. Under good conditions on $\Cc$, each object $E$
 gives a $\k$-point $[E]\in\Men_\Cc$ and we have the Kodaira-Spencer quasi-isomorphism
 \[
 T^\bullet_{[E]} \Men_\Cc \,\,\simeq \,\, R\Hom_\Cc(E,E)[1]. 
 \]
 
 \noindent {\bf ``Symplectic":} understood in the sense of \cite{PTVV}. That is, the datum of
 a symplectic form on a derived stack $\Yen$ includes not only pairings on the tangent dg-spaces
 $T^\bullet_y \Yen$ but also higher homotopies for the de Rham differentials of such pairings. 
 
 \vskip .2cm
 
 \noindent {\bf ``2-Calabi-Yau":} In order for the approach of \cite{kontsevich-soibelman} to be applicable,
 even at the formal level, 
 we need not only canonical identifications
 $R\Hom(E,F)^* \simeq R\Hom(F,E)[2]$ but a finer structure: a class in the Hochschild
 cohomology of $\Cc$ inducing these identifications. For instance, if $\Cc$ is
 the derived category of dg-modules over a dg-algebra $L$, we need an isomorphism
 $\gamma$ as in \eqref{eq:hoch}, i.e., $L$ should have a structure of a Calabi-Yau dg-algebra
 in the sense of \cite{Gi}. For the categories of deformation quantization 
 modules, Hochschild cohomology classes of this 
 nature were constructed in \cite{KS-DQ} Thm. 6.3.1.
 
 \vskip .2cm
 
  While there is every reason to expect the validity of  Principle \ref{CYP}
  in this setting, this has not yet been established.
  The case of $\Cc = D\Mc(X,\emptyset) = D^b_{\text {loc. const}}(X)$ for a smooth compact
  $X$ follows from the results of \cite{PTVV}, as in this case $\Men_\Cc$
  is interpreted in terms of   mapping stacks to the $(-2)$-shifted symplectic stacks $BGL_N$.
  This interpretation does not apply to $D\Mc(X,\emptyset)$ for a general
  compact nodal curve $X$. 
 So we cannot use Principle \ref{CYP}
  to construct ``symplectic moduli spaces of microlocal sheaves". In the next section
  we present an alternative, more direct approach via quasi-Hamiltonian reduction.

   \vfill\eject
   
   \section{Framed microlocal sheaves and multiplicative quiver varieties}\label{sec:framed}
   
   \sparagraph {\bf Motivation.} Recall \cite{N} that the original setting
   of Nakajima Quiver Varieties $M_\Gamma(V,W)$
   involves two types of vector spaces
   associated to vertices $i$ of  quiver $\Gamma$:
   \begin{enumerate}
   \item[(1)] The ``color" spaces $V_i$ which are {\em gauged},
   i.e., we perform the Hamiltonian reduction by the group
    $GL(V)=\prod GL(V_i)$ in order to arrive at $M_\Gamma(V,W)$.
   
   \item[(2)] The ``flavor" spaces $W_i$ which are 
   {\em fixed}, in the sense that $M_\Gamma(V,W)$ depends on $W$ functorially.
   In particular, it has a Hamiltonian action of the
    group $GL(W)=\prod GL(W_i)$.
    
    The setting of preprojective algebras (whose multiplicative version
was reviewed in \S \ref{sec:prepro}), corresponds to the case when
 $W_i=0$.

   \end{enumerate}

   \noindent In this section we explain a geometric framework allowing us to introduce
   such flavor spaces in a more general context of microlocal sheaves.
   For simplicity we restrict   the discussion to
   untwisted microlocal sheaves.  
   
   \vskip .2cm
   
   \sparagraph  {\bf Microlocal sheaves framed at $\infty$.}
   Let $Y$ be a quasi-projective nodal curve over $\CC$ with a duality 
   structure. We assume that 
   $Y=\ol {Y} -\infty$, where $\ol Y$ is a compact nodal curve  and
  $\infty =  \{\infty_j\}_{j\in J}$ is a finite set of smooth points. 
  Let
   \[
   Y^\partial \,\,=\,\,\on{Bl}_{\infty}(\ol Y) \,\,=\,\, Y \sqcup C, 
   \quad C \,\,=\,\, \bigsqcup_{j\in J} C_j
   \]
   be the real blowup of $\ol Y$ at $\infty$. Thus $Y^\partial$ is a compact topological
   space obtained by adding to $Y$ the circles $C_j$, so that each $C_j=S^1_{\infty_j}\ol Y$
    is the circle
   of real directions of $\ol Y$ at $\infty_j$. Note that in a neighborhood of $C$, 
    the space  $Y^\partial$  is naturally a 2-dimensional oriented
   $C^\infty$-manifold with boundary $C$. 
    We choose a base point $p_j$ in each $C_j$.
    
   Any microlocal sheaf  $\Fc$ on $Y$ is a local system in degree 0 near $\infty$.
   Thus it extends canonically 
  (by direct image) to  a complex of sheaves $\Fc^\partial$ on $Y^\partial$
   which is a local system in degree 0 near $C$. In particular, it gives rise to
 finite-dimensional  $\k$-vector spaces $\Fc_{p_j}$,  defined as the stalks of $\Fc^\partial$ at $p_j$. 
     We denote by
   \[
   \men_j(\Fc): \Fc_{p_j}\lra\Fc_{p_j}
   \]
   the  anti-clockwise monodromy of $\Fc^\partial$ around $C_j$

   \begin{Defi}\label{def:framed}
   Let $W=(W_j)_{j\in J}$ be a family of finite-dimensional $\k$-vector spaces.
   By a $W$-{\em framed microlocal sheaf} on $Y$ we mean a datum consisting
   of a microlocal sheaf $\Fc\in\Mc(Y,\emptyset)$  together with
   isomorphisms $\phi_j: \Fc_{p_j} \to W_j$, $j\in J$.

   We denote by $\Mc(Y)_W$ the category (groupoid) formed by $W$-framed
   microlocal sheaves on $Y$ and their isomorphisms (identical on $W$).
    
  \end{Defi}

     \begin{prop}\label{prop:affine-mod}
   Assume that $Y$ is an affine nodal curve with a duality structure,
    i.e., there is at least one
   puncture on each irreducible component. Then:
   \begin{itemize}
  \item[(a)] There exists a smooth affine algebraic $\k$-variety $\Men(Y)_W$
  (the {\em moduli space of $W$-framed microlocal sheaves})
   such that  isomorphism
  classes of objects of $\Mc(Y)_W$  are in bijection with
  $\k$-points of $\Men(Y)_W$.
  
  \item[(b)] The group $GL(W)=\prod GL(W_j)$ acts on $\Men(Y)_W$ by change of
   the framing. Taking the monodromies around the $C_j$  
  gives  an equivariant
   morphism (which we call the {\em moment map})
   \[
   \men = (\men_j)_{j\in J}: \Men(Y)_W \lra GL(W). 
   \]
  \end {itemize}   
   
  \end{prop}
  
  \noindent {\sl Proof:} (a)
We analyze the data of a $W$-framed microlocal sheaf  directly
  on $\wt X$, as in the previous section. These data reduce to  
  a collection of linear operators between the $W_j$ such that
  certain expressions formed out of them are invertible but, since
  each $\wt X_i$ is affine, subject to no other relations.
  This means that $\Men(Y)_W$ is realized as an open subset
  in the product of sufficiently
   many copies of affine spaces $\Hom(W_j, W_{j'})$.
  
  (b) Obvious. 
  
  \qed

   \begin{ex}[(Smooth Riemann surface)]\label{ex:RS}
  (a)  Let $\ol Y$ be a smooth projective curve of genus $g$.
   Choose one point $\infty\in\ol Y$ and put $Y=\ol Y-\{\infty\}$,
   so that $|J|=1$. Accordingly, we choose one base point $p\in Y$
   near $\infty$ in the sense explained above. A microlocal sheaf
   $\Fc\in\Mc(Y,\emptyset)$ is just a local system on $Y$. 
   
   So we fix one vector space $W$ and denote $G=GL(W)$. 
   A $W$-framed microlocal sheaf is just a homomorphism
   $\pi_1(Y, p)\to G$. As well known, $\pi_1(Y,p)$ is a free group
   on $2g$ generators $\alpha_1, \cdots, \alpha_g, \beta_1, \cdots, \beta_g$
   which correspond to the $a$- and $b$-cycles on the compact curve $\ol Y$. 
   Therefore $\Men(Y)_W = G^{2g}$ is the product of $2g$ copies of $g$.
    The $G$-action on $\Men(Y)_W$
   is by simultaneous conjugation. The moment map has the form
   \[
   \men: G^{2g} \lra G, \quad (A_1, \cdots,  A_g, B_1, \cdots,  B_g)
   \,\mapsto \, \prod_{\nu=1}^g [A_\nu, B_\nu],
   \]
   so $\men^{-1}(e) = \Hom(\pi_1(\ol Y, \infty), G)$ is the set of local systems
   on the compactified curve, trivialized at $\infty$. 
   
   \vskip .2cm
   
   (b) More generally, let $Y$ be an arbitrary smooth curve, 
   compactified to $\ol Y$ by a finite set of punctures
   $\infty_j, j\in J$. Then $\Men(Y)_W$
   is the space of representations of $\pi_1(Y, \{\infty_j\}_{j\in J})$,
   the  fundamental groupoid of
   $Y$ with respect to the set of   base points $\infty_j$. This is
   the setting of \cite {AMM},  \S 9.2, see also \cite{Boa}, Thm. 2.5. 
   
   \end{ex}
   
   \begin{ex}[(Coordinate cross)]\label{ex:cross}
    Let $Y=\{(x_1,x_2)\in\AAA^2 | \,\, x_1 x_2=0\}$ be the union of two affine lines meeting transversely.
    Then $\ol Y$ is the union of two projective lines meeting transversely and $\infty$ consists of two
    punctures. Accordingly, we have two marked points on $Y^\partial$, denote them $p_1$ and $p_2$.
 Given a family of two vector spaces $W=(W_1, W_2)$, the stack $\Men(Y)_W$
 is the  affine algebraic variety known as the 
{\em  van den Bergh's quasi-Hamiltonian space}, see \cite{vdB} and \cite[\S 2.4]{Boa}:
 \[
 \Men(Y)_W \,\,=\,\,  \Bc(W_1, W_2) \,\, := \,\, 
  \bigl\{ \xymatrix{
 W_1
   \ar@<.5ex>[r]^{a}&  \ar@<.5ex>[l]^{b}
  W_2
  } \bigl| 1+ab \text{ is invertible}\bigr\}.
 \]
  Note that $1+ba$ is also invertible on $\Bc(W_1, W_2)$. 
   \end{ex}

   \begin{ex}[(Microlocal sheaves with framed $\Phi$)]\label{ex:fr-phi}
    Let $X$ be a  compact nodal curve  with a duality structure,
   and $A\subset X$ be a finite
   subset of smooth points.
    Form a new curve $Y=X_A$, as in
   Proposition \ref{prop:A>empty}. Then $\Mc(Y)_W$ can be seen as the
   category parametrizing 
   microlocal sheaves on $X$ which are allowed singularities at $A$,
   but are equipped with a $W$-framing of their vanishing cycles
   at each such singular point. To emphasize it, we denote this category
   by $\Mc(X,A)_W$. 
  
    \end{ex}

    \begin{ex}[(Multiplicative quiver varieties)]
    We now specialize the above example further. 
    Let $X$ be a compact nodal curve with irreducible components
    $X_i, i\in I$. Assume, as in \S \ref{sec:prepro}, that each $X_i$
    is a rational curve, i.e., that the normalization $\wt X_i$
    is isomorphic to $\PP^1$. Choose the set $A$ consisting
    of precisely one smooth point $a_i$ on each $X_i$. 
    Let $W=(W_i)_{i\in I}$ be a family of $\k$-vector spaces. 
    Thus the topological structure of $(X,A)$ is determined by the graph $\Gamma$
    of intersections of irreducible components of $X$, in particular, $I$ is the set of vertices of $\Gamma$. 
    We will write $X=X_\Gamma$
    to indicate this dependence.
    
    \end{ex}

   \begin{prop}
   In the situation just described,  $\Mc(X,A)_W$
   is  equivalent to the category which parametrizes  linear algebra
   data consisting of:

   \begin{enumerate}
   \item[(1)] Collections of vector spaces $V=(V_i)_{i\in I}$; 
   
   \item[(2)] Linear maps 
   \[
   \begin{gathered}
   a_h: V_{s(h)}\to V_{t(h)}, \quad b_h: V_{t(h)}\to V_{s(h)}, \quad h\in E,
   \\
   u_i: V_i\to W_i, \quad v_i: W_i\to V_i, \quad i\in I, 
  \end{gathered}
   \]
    such that all the maps 
   \[
   (\1+a_h b_h),\,\, (\1+b_h a_h), \,\, (\1+u_iv_i), \,\, (\1+v_iu_i)
   \]
   are invertible, and
   \item[(3)]
   For each $i\in I$ we have the identity
    \[
(\1_{V_i}  + v_i u_i)  \prod_{h\in E, \, t(h)=i} (\1_{V_i} +
 a_h b_h) \prod_{h\in E, \, s(h)=i} (\1_{V_i} + b_h a_h)^{-1} 
  \,\,=\,\, \1_{V_i}.
  \]
 \end{enumerate}
  These data are considered modulo isomoprhisms of the $V_i$.

   \end{prop}

   \noindent{\sl  Proof:}  Completely analogous to that of Theorem \ref{thm:prepro}
   and we leave it to the reader. \qed

   \vskip .2cm
   
     The moduli spaces of semistable objects of $\Mc(X,A)_W$ (\
   defined as GIT quotients) as well as their analogs for twisted sheaves are the
  {\em multiplicative quiver varieties} (MQV) as defined in  \cite{yamakawa}.

   \begin{ex}[(Higher genus MQV)]\label{ex:hi-MQV}
   In the interpretation of the previous example, we associated  to a graph $\Gamma$  a  nodal curve
 $X_\Gamma$ with all components rational. In particular, the number $g_i$ of loops at a vertex $i\in\Gamma$
 was interpreted as the number of self-intersection points of the corresponding rational curve $X_i$. 
 We can also associate  to $\Gamma$ a nodal curve $X'_\Gamma$ in a different way,
 by taking the component $X'_i$ associated to $i$ to be of genus $g_i$
 (and interpreting other edges of $\Gamma$ as intersection points of the $X'_i$). 
 Choose the set $A$ to consist of one point on each irreducible component of $X'_\Gamma$. 
 This defines a datum $(X'_\Gamma, A)$ uniquely up to a diffeomorphism. 
 We will refer to the moduli spaces of objects of $\Mc(X'_\Gamma, A)_W$
 (defined as GIT quotients) as {\em higher genus multiplicative quiver varieties}
 associated to $\Gamma$. Note that one can also consider their twisted versions, involving
 twisted microlocal sheaves.   
   
   \end{ex}

 \sparagraph {\bf Quasi-Hamiltonian $G$-spaces.} Here we review the main points of the
 theory of group valued moment maps  \cite{AMM}.   For simplicity we work in
 the complex algebraic situation, not that of compact Lie groups.
 
 Let $G$ be a reductive algebraic group over $\CC$, with Lie algebra $\gen$. 
We denote by
\[
\theta^L = g^{-1}dg, \,\, \theta^R = (dg)g^{-1} \,\,\in \,\,\Omega^1(G, \gen)
\] 
 the standard left and right invariant $\gen$-valued 1-forms on $G$. 
 
 We fix an invariant symmetric bilinear form $(-,-)$ on $\gen$. It gives rise to
 the bi-invariant scalar 3-form (the {\em Cartan form})
 \[
 \eta = {1\over 12} (\theta^L, [\theta^L, \theta^L]) \,\,=\,\, 
  {1\over 12} (\theta^R, [\theta^R, \theta^R])\,\,\in \,\,\Omega^3(G).
 \]
 For a $G$-manifold $M$ and $\xi\in\gen$ we denote by $\partial_\xi$ the
 vector field on $M$ corresponding to $\xi$ by the infinitesimal action.

 \begin{Defi} \label{def:AMM}  \cite{AMM}  
 A quasi-Hamiltonian $G$-space is a smooth algebraic variety  $M$ with $G$-action,
  together with a $G$-invariant
 2-form $\omega\in\Omega^2(M)^G$ and a $G$-equivariant map
 $\men: M\to G$ (the {\em group valued moment map}) such that:
 \begin{enumerate}
 \item[(QH1)] The differential of $\omega$  satisfies  $d\omega = -\men^*\chi$. 
 
 \item[(QH2)] The map $\men$ satisfies, for each $\xi\in\gen$, the condition
 \[
 i_{\partial_\xi} \omega \,\,= {1\over 2} \men^*(\theta^L+ \theta^R, \xi). 
 \]
 Here $(\theta^L+ \theta^R, \xi)$ is the scalar 1-form on $G$ obtained by
 pairing the $\gen$-valued form $\theta^L+ \theta^R$ and the element $\xi\in\gen$
 via the scalar product  $(-,-)$. 
 
 \item[(QH3)] For each $x\in M$, the kernel of the 2-form $\omega_x$ on $T_xM$ is
 given by
 \[
 \on{Ker}(\omega_x) \,\,=\,\, \bigl\{ \partial_\xi (x),
  \,\,\xi\in \on{Ker}(\on{Ad}_{\men(x)} +\1\bigr\}. 
 \]
 \end{enumerate}
 
 \end{Defi}

Given a quasi-Hamiltonian $G$-space $(M,\omega, \men)$,  the
 {\em quasi-Hamiltonian reduction} of $M$ by $G$ is, classically \cite{AMM},   
 the   quotient
 \[
  M\qqq G \,\,=\,\, \men^{-1}(e)^{\on{sm}} / G, 
 \]
 where $\men^{-1}(e)^{\on{sm}}$ is  the smooth locus of the scheme-theoretic
 preimage $\men^{-1}(e)$ or, more precisely, the open part formed by those points $m$, 
 for which $d_m\men$ is surjective.   
 \begin{thm}\label{thm:qred} \cite{AMM}
 For any quasi-Hamiltonian $G$-space $M$ the quotient $M\qqq  G$ is a smooth
 orbifold (i.e., Deligne-Mumford stack) with 
a canonical symplectic structure. \qed
 \end{thm}

 \begin{rem}
 Using the language of derived stacks allows one to formulate Theorem \ref{thm:qred}
 in a more flexible way, without restricting to the locus of smooth points. More precisely,
 we can form the smooth derived stack of amplitude $[-1,1]$
 \[
 [M\qqq G]^{\on{der}}  \,\,=\,\, R\men^{-1}(e) \qq G. 
 \] 
 Here $R\men^{-1}(e)$ is the derived preimage of $e$, a smooth derived scheme
 of amplitude $[0,1]$.  Further, the symbol $-\qq G$ means stacky quotient by $G$. The analog of
 Theorem \ref{thm:qred} is then that $[M\qqq G]^{\on{der}}$ is a symplectic derived
 stack which contains $M\qqq G$ as an open part. 
 \end{rem}
 
 The following is the main result of this section. 
 
 \begin{thm}\label{thm:MW-quasi}
 Let $Y$ be an affine nodal curve, and $W=(W_j)$ as before. The smooth algebraic
 variety $\Men(Y)_W$ has a natural structure of a quasi-Hamiltonian $GL(W)$-space
 with the moment map $\men=\men_W$
  given by the monodromies (Proposition \ref{prop:affine-mod}(b)). 
  \end{thm}
  
 \begin{rem} This result provides a more direct approach to the ``moduli space"
 of microlocal sheaves on a compact nodal curve, in particular, 
 to the symplectic structure on this space. 

Indeed, the set-theoretic quotient $\men_W^{-1}(e)/GL(W)$ parametrizes microlocal
sheaves $\Fc$ on the compact curve $\ol Y$ such that  the dimensions of the stalk
 of $\Fc$ at  $\infty_j$ is equal to  $\dim W_j$. Thus we can {\em define}
 the derived stack
 \[
 \Men(\ol Y, \emptyset) \,\,=\,\, \bigsqcup_{W}\,\,  [\Men(Y)_W\qqq GL(W)]^{\on{der}},
 \] 
 the disjoint union over all possibe choices of $(\dim W_j)_{j\in J}$. 
 
 Alternatively, one can consider the Poisson variety obtained as the spectrum
 of the algebra of $GL(W)$-invariant functions on $\Men(Y)_W$, cf. \cite{Boa},
 Prop. 2.8. 
 
 \end{rem}

\vskip .2cm
In the case of a smooth curve $Y$, see   Example \ref{ex:RS}(b), a proof of
Theorem 
   \ref{thm:MW-quasi}  
was given in   \cite[\S 9.3] {AMM} using a procedure called {\em fusion} which
allows one to construct complicated quasi-Hamiltonian spaces from simpler ones.
We use the same strategy but allow one more type of ``building block" in the fusion
construction.

 \vskip .2cm
 
 \sparagraph {\bf  Fusion of quasi-Hamiltonian spaces.}
  We now briefly review the  necessary concepts. 
 
 \begin{thm}[\cite{AMM}]\label{thm:fusion}
 Let $M$ be a quasi-Hamiltonian $G\times G\times H$-space, with moment map
$\men=(\men_1, \men_2, \men_3)$. Let $G\times H$ act on $M$ via the
diagonal embedding $(g,h) \mapsto (g,g,h)$. Then $M$ with the 2-form
\[
 \omega' = \omega + {1\over 2} (\men_1^*\theta^L, \men_2^*\theta^R)
\]
and the moment map
\[
\men' \,\,=\,\,  (\men_1\cdot \men_2, \men_3): M\lra G\to H 
\]
is a quasi-Hamiltonian $G\times H$-space, called the (intrinsic) {\em fusion}
of the $G\times G\times H$-space $M$. 
\end{thm}

 \begin{rem}
The geometric meaning of the fusion is that the two copies of $G$
from $G\times G\times H$  are ``attached"
to the two of the tree boundary components of a 3-holed sphere, and the new
diagonal  copy of $G$ from $G\times H$ is then ``read off" on the remaining component,
see \cite{AMM}, Ex. 9.2 and  \cite{Boa} \S 2.2. Thus, in the case of smooth curves, fusion directly
corresponds to gluing Riemann surfaces out of simple pieces. We will extend this
to nodal curves. 
 \end{rem}
 
 The {\em extrinsic fusion} of a quasi-Hamiltonian $G\times H_1$-space $M_1$ and 
 a $G\times H_2$-space $M_2$ is the $G\times H_1\times H_2$-space
$M_1 \circledast M_2$ which is the fusion of the $G\times H_1\times G\times H_2$-space
$M_1\times M_2$ along the embedding $G\to G\times G$. 

We will use the following three building blocks.

\begin{exas}
(a) {\bf (Double of $G$: annulus).} Given $G$ as before, its {\em double}
is the quasi-Hamiltonian $G\times G$-space $D(G)=G\times G$ with coordinates $a,b\in G$,
the $G\times G$-action given by
\[
(g_1, g_2) (a,b) = (g_1 a g_2^{-1}, g_2 b g_1^{-1}),
\] 
the moment map given by
\[
\men_D: D(G) = G\times G \lra G\times G, \quad (a,b) \mapsto (ab, a^{-1}, b^{-1})
\]
and the 2-form given by
\[
\omega_D = {1\over 2} (a^*\theta^L, b^*\theta^R) + {1\over 2} (a^*\theta^R, b^*\theta^L).
\]
For a vector space $V$ and  $G=GL(V)$, this space is identified with
$\Men(Y)_W$, where $Y$ is a 2-punctured sphere and $W=(V,V)$ associates $V$ to each puncture.
The surface with boundary $Y^\partial$ is an annulus. 

\vskip .2cm

(b) {\bf Intrinsically fused double: holed torus.} With $G$ as before, its
{\em intrinsically fused double} ${\bf D}(G)$ is the quasi-Hamiltonian $G$-space $G\times G$
obtained as the fusion of the $G\times G$-space $D(G)$. For a vector space $V$ and
$G=GL(V)$, this space is identified with $\Men(Y)_V$ where $Y$ is a 1-punctured elliptic curve.
The surface $Y^\partial$ is  a 1-holed torus.

\vskip .2cm

(c) {\bf The space $\Bc(W_1,W_2)$: cross.} To treat nodal curves, we add the third type of
 building blocks: the varieties $\Bc(W_1,W_2)$, see Example
 \ref{ex:cross}. Again, this is a known quasi-Hamiltonian $GL(W_2)\times GL(W_2)$-space
  \cite{vdB} \cite{vdB2} with moment map
  \[
  (a,b) \,\,\mapsto\,\, \bigl((1+ab)^{-1}, 1+ba\bigr) \,\,\in 
  \,\, GL(W_2)\times GL(W_1)
  \]
  and the 2-form
  \[
  \omega \,\,=\,\, {1\over 2} \bigl( \on{tr}_{W_2}(1+ab)^{-1}da \wedge db-
   \on{tr}_{W_1}(1+ba)^{-1} db\wedge da\bigr). 
  \]
  As we saw in Example \ref{ex:cross}, it is identified with $\Men(Y)_{W_1, W_2}$,
  where $Y$ is a coordinate cross. The  topological space  $Y^\partial$ is 
  the union of two disks meeting at one point. 
  
\end{exas}

Let now $Y$ be an affine nodal curve.
The topological space $Y^\partial$ can then
be decomposed into elementary pieces of types (a)-(c) in the above examples,
joined together by several 3-holed spheres.

Let  $W=(W_j)_{j\in J}$ be given.
Note that for $\Men(Y)_W$ to be non-empty,  the numbers $N_j = \dim W_j$ should depend only on
the irreducible component of $Y$ containing $\infty_j$.
This means that to each boundary component of each elementary piece
we can unambiguously associate a group $GL(N_j)$ and so form the corresponding
quasi-Hamiltonian space of type (a)-(c) above. Taking the product of
these corresponding quasi-Hamiltonian spaces and performing the fusion along
the 3-holed spheres, we get a quasi-Hamiltonian space which is identified with
$\Men(Y)_W$. This proves Theorem \ref{thm:MW-quasi}.

\begin{rem}
It would be interesting to construct the 2-form on $\Men(Y)_W$ more intrinsically,
in terms of a cohomological pairing, using  some version of Poincar\'e-Verdier
duality for cohomology with support on the ``nodal surface with boundary" $Y^\partial$. 
This does not seem to be known even  for  smooth $Y$. 
\end{rem}

\vfill\eject

\section{Further directions}

\sparagraph{(Geometric) Langlands correspondence for nodal curves.}
Since microlocal sheaves without singularities form a natural analog
of local systems for nodal curves, it would be interesting to put them
into the framework of the Langlands correspondence. In particular,
for a compact nodal curve $X$
it would be interesting to have a  derived equivalence
between the  de Rham version (cf. \S \ref{sec:dr})
of the  ``Betti-style" derived  stack $\Men(X,\emptyset)$ and  some other moduli stack
$\mathfrak B$ 
of ``coherent" nature, generalizing the moduli stack of vector bundles
on  a smooth curve. A potential candidate for $\mathfrak B$  is provided by the moduli stack
of Riemann surface quiver representations in the sense of \cite{CB}. 

Note that the concept of microlocal sheaves makes sense for nodal curves $X$ over
$\FF_q$.  So one can expect that their $L$-functions
(appropriately defined)
have, for projective nodal curves $X$, properties similar to those of
$L$-functions of local systems on smooth projective curves over $\FF_q$.

One can even consider arithmetic analogs of nodal curves, obtained
by  gluing the spectra of rings of integers in number fields
along closed points. An example is provided by the spectrum of the
group ring $\ZZ[\ZZ/p]$, where $p$ is a prime. This scheme is the union
of $\Spec(\ZZ)$ and $\Spec (\ZZ[\sqrt[p] {1} ])$ meeting transversely at the point $(p)$,
cf. \cite{Mi}, \S 2.

\sparagraph {Multiplicative quiver varieties and mirror symmetry.}
 Let $\Gamma$ be a finite graph, possibly with loops, and 
  $\mathbb{M}_{V,W}(X'_\Gamma)^q$ be the corresponding
   {\em  higher genus multiplicative quiver varieties}, see
 Example \ref{ex:hi-MQV}. Here $q=(q_i)\in(\CC^*)^I$ is a vector of twisting parameters. 
 We expect that the varieties   $\mathbb{M}_{V,W}(X'_\Gamma)^q$ are mirror dual to 
the ordinary (``additive")  Nakajima quiver varieties for $\Gamma$. 

In particular, we expect  that   $\mathbb{M} _\Gamma(V,W)^q$
is singular  if and only if the point $q$
lies in the singular locus of the equivariant quantum connection
for the ordinary quiver variety.
Here,  equivariance is in reference to the action of an algebraic torus which acts on the quiver
variety scaling the symplectic form by a nontrivial character. See    \cite{MO},  where this
  connection  as well as its singularities,  have been computed. 

\sparagraph {Borel/unipotent reduction and cluster varieties.}
It would be interesting to extend the approach of \cite{FG}
from local systems on smooth curves to microlocal sheaves
on nodal curves. 
  That is, 
in the situation
of \S \ref{sec:framed}B we can  choose any number of marked points $p_{j,\nu}$ on each boundary
component $C_j$ of $Y^\partial$. After this we  can consider microlocal sheaves $\Fc$
 together with a Borel or unipotent reduction of the structure group
 at each $p_{j,\nu}$ (recall that each restriction $\Fc|_{C_j}$ is a local system). 

 This can lead to interesting cluster varieties. These varieties may be related
 to the  classification of irregular DQ-modules on  a symplectic surface with
 support in a nodal curve.

\sparagraph{3-dimensional generalization.}  
The datum of a smooth compact curve over $\CC$ (topologically, an 
oriented surface) $X$
and a finite set of points $A\subset X$ has the following 3-dimensional analog.

We consider a compact oriented $C^\infty$ 3-manifold $M$ and a {\em link}
in $M$, i.e., a collection $L=\{C_a\}_{a\in A}$ of disjoint embedded
circles ({\em knots}).
 We have then a stratification of $M$ into the $C_a$ and the complement
of their union. Denote the $D^b_L(M)$ the category of complexes 
of sheaves on $M$, constructible
with respect to this stratification. For  $L=\emptyset$,  it is a 3-Calabi-Yau category by Poincar\'e duality. 
For  arbitrary $L$, it has a natural abelian
subcategory $\Perv(M, L)$  of ``perverse sheaves". Given any surface $X\subset M$
meeting $L$ transversely,
an object $\Fc\in\Perv(M,L)$  gives a perverse sheaf on $X$, smooth outside 
  $X\cap L$. 

One can obtain 3d analogs of compact nodal curves (``nodal 3-manifolds")
 by identifying several compact 
3-manifolds pairwise along some knots.  For example, we can glue two
such manifolds $M'$ and $M''$ (say, two copies of the sphere $S^3$)
 by identfying a knot $C'\subset M'$
with a knot $C''\subset M''$. As the normal bundle $T_C M$  of a knot 
$C$ in an oriented
3-manifold $M$ is trivial, we can choose  a duality structure,
i.e.,  an identification
of $T_{C'}M'$ with $T^*_{C''} M''$, and then set up the formalism of
microlocal  sheaves  and complexes. This should  lead to interesting 3-Calabi-Yau categories
and to $(-1)$-shifted symplectic stacks parametrizing their objects. 
  
   3-Calabi-Yau categories of the form $D\Mc^\Lc(X,\emptyset)$, see
   Theorem \ref{thm:CY-twist}(a), 
    correspond to a particular type of  nodal 3-manifolds:
    circle bundles over  nodal curves over $\CC$.

  \vfill\eject
  
  \appendix

\section { Notations and conventions.} 
 We fix a base field $\k$. All sheaves will be understood
 as   sheaves of $\k$-vector spaces, similarly for complexes of sheaves. 
 
 All topological spaces we consider will be understood to be homeomorphic to
 open sets in finite CW-complexes,
 in particular,  they are locally compact and of finite dimension.
 For a space $X$ we denote by $\Sh(X)$ the category of sheaves of $\k$-vector spaces
 on $X$. We denote by $\DSh(X)$ the bounded derived category of $\Sh(X)$.
 We will consider it as a pre-triangulated category \cite{BK},  i.e.,   as a dg-category enriched by the complexes $R\Hom(\Fc, \Gc)$, so that $H^0R\Hom(\Fc,\Gc)$
 is the ``usual" space of morphisms from $\Fc$ to $\Gc$ in the derived category. 
 Alternatively, we can view it as a stable $\infty$--category by passing to the dg-nerve  \cite{lurie-stable} \cite{lurie-algebra} \cite{faonte}. 
 
 We denote by $\Dbcc(X)\subset \DSh(X)$ the full subcategory of cohomologically constructible complexes \cite{KS} and by $\DD=\DD_X$
 the Verdier duality functor on this subcategory   \cite[\S 3.4]{KS}. Thus, if $X$ is an oriented $C^\infty$-manifold of real dimension $d$, and 
 $\Fc$ is a local system on $X$ (put in degree 0), then $\DD(\Fc) = \Fc^\bigstar[d]$, where $\Fc^\bigstar$ is the dual local system.
 In general, for any compact space $X$ and any $\Fc\in\Dbcc(X)$ we have  {\em Poincar\'e-Verdier duality},
which is the canonical  identification of complexes of $\k$-vector spaces with finite-dimensional
cohomology, and consequently, of their cohomology spaces:
  \be\label{eq:PVD}
  \begin{gathered}
  R\Gamma(X,\Fc)^*  \,\, \simeq\,\,  R\Gamma(X, \DD_X(\Fc));
  \\
 \HH^i(X, \Fc)^* \,\, \simeq\,\,  \HH^{-i}(X, \DD_X(\Fc)).
 \end{gathered}
 \ee
 
 Let   $X$  be a complex manifold. We denote
 by $\Dbc(X)\subset\Dbcc(X)$ the derived
 category of bounded complexes of sheaves on $X$ with $\CC$-constructible cohomology sheaves. 
 The functor $\DD_X$ preserves this subcategory. 
 We denote by $\Perv(X)\subset\Dbc(X)$ the subcategory of perverse sheaves. 
   The conditions of perversity are normalized
 so that a local system on $X$  put in degree 0, is perverse. Thus $\Perv(X)$
 has the perfect duality given by
 \[
 \Fc\mapsto \Fc^\bigstar \,\, := \,\, \DD(\Fc)[-2\on{dim}_\CC(X)]. 
 \]

 \vfill\eject

\vskip .4cm

R. B.: Department of Mathematics, MIT, Cambridge MA 02139 USA,
{\tt bezrukav@math.mit.edu}

\vskip .2cm

M.K.:  Kavli IPMU, University of Tokyo, 5-1-5 Kashiwanoha, Kashiwa, Chiba,
277-8583 Japan,  {\tt mikhail.kapranov@ipmu.jp }

\end{document}